\documentclass{siamart190516}
\usepackage{amsfonts}
\usepackage{epstopdf}
\usepackage{algorithmic}
\usepackage{verbatim}
\usepackage[all]{xy}
\usepackage{graphicx}
\usepackage{float}
\CompileMatrices

\usepackage[shortlabels] {enumitem}
\usepackage{xcolor}
\usepackage{hyperref}
\hypersetup{
    colorlinks   = true,
    citecolor    = blue,
    urlcolor=blue,
}

\usepackage{amsopn}

\newcommand{\R}{{\mathbb R}}

\newcommand{\abs}[1]{\left| #1 \right|}

\newcommand{\norm}[1]{\left|\left|#1\right|\right|}

\newcommand{\setint}{{\mathrm{int}\,}}

\newcommand{\argmin}{\mathop{\mathrm{argmin }}}

\newcommand{\prox}{\mathop{\mathrm{prox}}}
\newcommand{\diag}{\mathop{\mathrm{diag}}}
\newcommand{\sign}{\mathop{\mathrm{sign}}}
\newcommand{\Var}{\mathop{\mathrm{Var}}}
\newcommand{\rank}{\mathop{rank}}

\usepackage{mathabx}
\usepackage{subfig}
\newcommand{\proxi}[2]{\prox_{\quad\quad#1}^{\quad\quad#2}}

\title{One-Step Estimation with Scaled Proximal Methods\footnote{\textbf{Funding: }Bassett was supported in part by ONR grants N0001419WX00183 and N0001420WX01523. Deride was supported by CONICYT-Chile through FONDECYT grant number 11190549.}}
\author{Robert Bassett\thanks{Dept. of Operations Research, Naval Postgraduate School, CA
(\email{robert.bassett@nps.edu}).}
\and Julio Deride\thanks{Department of Mathematics, Universidad Federico Santa Mar\'ia (\email{julio.deride@usm.cl}).}}

\headers{One-Step Estimation with Scaled Proximal Methods}{Robert Bassett
and Julio Deride}
\begin{document}
\maketitle
\begin{abstract}
We study statistical estimators computed using iterative optimization methods that are not run until completion. Classical results on maximum likelihood estimators (MLEs) assert that a \emph{one-step estimator} (OSE), in which a single Newton-Raphson iteration is performed from a starting point with certain properties, is asymptotically equivalent to the MLE. We further develop these early-stopping results by deriving properties of one-step estimators defined by a single iteration of scaled proximal methods. Our main results show the asymptotic equivalence of the likelihood-based estimator and various one-step estimators defined by scaled proximal methods. By interpreting OSEs as the last of a sequence of iterates, our results provide insight on scaling numerical tolerance with sample size. Our setting contains scaled proximal gradient descent applied to certain composite models as a special case, making our results applicable to many problems of practical interest. Additionally, our results provide support for the utility of the scaled Moreau envelope as a statistical smoother by interpreting scaled proximal descent as a quasi-Newton method applied to the scaled Moreau envelope.
\end{abstract}

%
%\section*{Todo}
%
%\begin{enumerate}
%\item Present a cohesive narrative of the local quadratic expansion/differentiability in quadratic mean.
%\begin{enumerate}
%\item Discretization trick
%\item Differentiability in quadratic mean. \cite{vanderVaart}[Theorem 5.39 and Theorem 5.23] seem relevant.
%\item Cite Ibragomovic and van der Vaart's results on nonsmooth model.
%\item Extension to resolvent
%\item Extension to $\hat{\theta}_{MLE}$ instead of $\theta_{0}$.
%\end{enumerate}
%\item Examples for unregularized version/proximal descent
%\begin{enumerate}
%\item Moreau Envelope as smoother
%\item Nonsmooth model and asymptotic efficiency
%\end{enumerate}
%\item Examples for regularized version/proximal gradient descent
%\begin{enumerate}
%\item Lasso/adaptive lasso/any estimator with oracle property \cite{zou}
%\item Generalized linear models \cite{zou}
%\item Matrix completion \cite{guihangu}[Corollary 3.6 and Remark 3.7] and \cite{negahbanwainwright}[Section 3.4.1].
%\end{enumerate}
%\end{enumerate}

\begin{keywords}
  Proximal Operator, One-Step Estimator, Moreau Envelope
\end{keywords}

% REQUIRED
\begin{AMS}
  62F12, 65K10, 90C53
\end{AMS}

\section{Introduction.} \label{Intro}

In likelihood-based statistical inference, estimators are defined as solutions to optimization problems, with the objective function constructed from a random sample. When computing the estimator, however, it is often the case that the statistical context for the optimization problem is disregarded and a purely numerical perspective adopted. Numerical tolerance is taken to be as small as can be computed in a reasonable amount of time, and the resulting approximate minimizer is taken as the realization of the estimator for the given data sample. 

In this paper we view statistical and numerical error holistically. We retain the statistical origin of the mathematical program defining the estimator in order to provide insight on the numerical tolerance required to achieve statistical optimality. By ``statistical optimality'' we mean asymptotic equivalence to the defined estimator, a minimizer which depends on a random sample. For reasons that will become clear in Section~\ref{Numerics}, we focus our attention on scaled methods. We consider scaled proximal gradient descent and scaled proximal descent, generalizations of proximal gradient descent and proximal descent, respectively, which adjust for curvature of the objective function. We show that one-step estimators constructed from these algorithms achieve statistical optimality in an asymptotic sense for a broad class of parametric families and likelihood-based estimators.

The rest of this paper is organized as follows. In Section~\ref{Intro}, we review one-step estimation for statistical inference in parametric models and summarize previous work on scaled proximal algorithms. In Section~\ref{Composite}, we give one-step estimation results for scaled proximal gradient descent applied to a composite model, such as is commonly encountered in penalized and constrained M-estimation. Section~\ref{ProxDesc} contains similar results for scaled proximal descent. The scaled and unscaled proximal operator can be interpreted as scaled and unscaled gradient descent, respectively, applied to an inf-convolution of the objective function, so in Section~\ref{Moreau} we provide an alternative interpretation of our results using variational analysis to describe the scaled Moreau envelope as statistical smoother. We conclude with Section~\ref{Numerics}, which provides a counterexample demonstrating that one-step estimation results of the type we consider do not hold for first order methods. In this section, we also provide numerical validation of our results and some examples demonstrating their utility.

\subsection{One-Step Estimation.} \label{subsec:one-step}

It is well known that for a general class of statistical models the MLE is asymptotically unbiased and efficient, meaning that as the number of samples goes to infinity, the expectation of the MLE matches the parameter to be estimated and its variance attains the Cramer-Rao lower bound. We begin by reviewing aspects of this theory relevant to our contributions.

Let $\Theta \subseteq \R^d$ be an open set and $\{P_{\theta}: \theta \in \Theta\}$ a parametric family of probability measures on a measurable space $(\mathcal{X}, \mathcal{A}, \mu)$, where $\mu$ is $\sigma$-finite. Assume that for all $\theta$ the measure $P_{\theta}$ is absolutely continuous with respect to $\mu$, hence has Radon-Nikodym derivative $p_{\theta} := dP_{\theta}/d\mu$, and its support does not depend on $\theta$. Fix $\theta^{*} \in \Theta$ and assume $X_{1},...,X_{n} \sim^{i.i.d.} p_{\theta^{*}}$. We will exclusively focus on estimators $\hat{\theta}$ which converge at the $\sqrt{n}$ parametric rate, so that $\sqrt{n}\left(\hat{\theta}-\theta^{*}\right)$ converges in law to some limit distribution. 

A sequence of random variables $X_n$ is stochasticially bounded, denoted $O_P(1)$, if $X_n$ is bounded in probability with respect to the measure $P_{\theta^{*}}$. Additionally, $X_n$ is $o_P(1)$ if it converges to zero in probability with respect to $P_{\theta^{*}}$. An estimator $\hat{\theta}$ of $\theta_{0}$ is said to be $\tau_{n}$-consistent, for a given sequence $\tau_{n}$, if $\tau_{n}(\hat{\theta} - \theta_{0})$ is $O_{P}(1)$. We also call $\hat{\theta}$ $\tau_{n}$-consistent for another estimator $\tilde{\theta}$ if $\tau_{n}(\hat{\theta} - \tilde{\theta}) = O_{P}(1)$. Because we only consider estimators which converge at the $\sqrt{n}$ parametric rate, two estimators $\hat{\theta}$ and $\tilde{\theta}_{2}$ are asymptotically equivalent when they are $\sqrt{n}$-consistent.

The statistical model $\{P_{\theta}: \theta \in \Theta\}$ is said to be \emph{differentiable in quadratic mean} at some fixed $\theta_{0} \in \Theta$ if there exists a measurable vector-valued function $\nabla{\ell}_{\theta_{0}}$ such that, as $\theta \to \theta_{0}$,
$$\int\left[ \sqrt{p_{\theta}} - \sqrt{p_{\theta_{0}}} - \frac{1}{2}(\theta - \theta_{0})^{\top} \nabla{\ell}_{\theta_{0}} \sqrt{p_{\theta_{0}}}\right]^{2} \, d\mu = o\left(\|\theta-\theta_{0}\|^{2}\right).$$

Differentiability in quadratic mean is a relaxed smoothness assumption which still permits many of the classical results related to maximum likelihood. For example, it can be shown that a location model with univariate Laplace density, where $p_{\theta} =\frac{1}{2} e^{-|x-\theta|}$ for $\theta \in \R$, is differentiable in quadratic mean with $\nabla{\ell}_{\theta}(x) = \sign(x-\theta)$ even though it is nonsmooth. For a statistical model that is differentiable in quadratic mean, the Fisher Information at $\theta_{0}$ is defined at $I_{\theta_{0}} = E\left[ \nabla{\ell}_{\theta_{0}} \nabla{\ell}_{\theta_{0}}^{\top} \right]$. The following theorem, due to Le Cam \cite{LeCam70}, gives the asymptotic properties of the maximum likelihood estimator; details can be found in \cite{vanderVaart}. %specifics for reference \cite{LeCam12}[Section 17.3], \cite{ibragimov}[Theorem and Lemma 1.1], \cite{vanderVaart}[Theorem 5.39].

\begin{theorem} \label{old_mle} %van de vaart thm 5.39
  Suppose that the model $\{P_{\theta}: \theta \in \Theta\}$ is differentiable in quadratic mean at $\theta^{*} \in \setint(\Theta)$. Suppose also that there exists a measurable function $L$ with $E[L^2(X_{1})]$ finite and such that, for every neighborhood of $\theta^{*}$ and every $\theta_{1}$ and $\theta_{2}$ in that neighborhood,
$$|\log p_{\theta_{1}}(x) - \log p_{\theta_{2}}(x)| \leq L(x) \|\theta_{1} - \theta_{2}\|,\quad\mbox{for every x in }\R^d$$
  If the Fisher information matrix $I_{\theta^{*}}$ is nonsingular and the MLE $\hat{\theta}_{mle}$ is consistent then $\sqrt{n}(\hat{\theta}_{mle} - \theta^{*})$ is asymptotically normal with mean zero and covariance matrix $I_{\theta^{*}}^{-1}$.
\end{theorem}

Theorem~\ref{old_mle} gives that the MLE is statistically optimal in a specific sense. According to the Cramer-Rao bound,
%\footnote{see for example theorem 7.3 in \cite{Ibragimov}},
 an unbiased estimator $\hat{\theta}$ of $\theta^{*}$ must satisfy the matrix inequality $\Var(\hat{\theta}) \succeq I_{\theta^{*}}^{-1}$. Because $\hat{\theta}_{mle}$ \emph{attains} this variance bound as $n \to \infty$, the MLE is optimal in terms of its asymptotic bias and variance. Moreover, this optimality of the MLE is the primary property which motivates its use. According to van der Vaart \cite{vanderVaart}, ``The justification through asymptotics appears to be the only general justification of the method of maximum likelihood.'' (See \cite{LeCam12} for a similar opinion). Though there are many important contributions related to finite-sample results for maximum likelihood estimation, \cite{FiniteWilks, Spokoiny}, they apply only with subgaussian or subexponential tail behavior. Our contributions are intended to apply to likelihood-based inference generally, so we will use asymptotic justifications for our results.

Asymptotic efficiency can be considered in the more general setting of M-estimation as well.  For each $\theta \in \Theta$, let $m_{\theta}: \mathcal{X} \to \mathbb{R}$ be a measurable function. The M-estimator for this criterion is defined as
$$\hat{\theta}_{M} = \argmin_{\theta \in \Theta} \frac{1}{n} \sum_{i=1}^{n} m_{\theta}(X_{i}).$$
If $m_{\theta} = -\log p_{\theta}$, the corresponding M-estimator coincides with the maximum likelihood estimator. M-estimation has an asymptotic result similar to Theorem \ref{old_mle} when $E\left[m_{\theta}(X)\right]$ has a second order expansion \cite[Theorem 5.23]{vanderVaart}.
%\begin{theorem}[\cite{vanderVaart}, Theorem 5.23] \label{VDV_M-est} For each $\theta \in \Theta$, let $m_{\theta}: \mathcal{X} \to \mathbb{R}$ be a measurable function such that $m_{\theta}(x)$ is differentiable at $\theta_{0}$ for $P$-almost every $x$
%%\footnote{We can also have $m_{\theta}$ differentiable at $\theta_{0}$ in $P$-probability.}
%with derivative $\nabla{m}_{\theta_{0}}(x)$ and such that, for every $\theta_{1}$ and $\theta_{2}$ in a neighborhood of $\theta_{0}$ and a measurable function $L$ with $P L^2 < \infty$
%$$|m_{\theta_{1}}(x) - m_{\theta_{2}}(x)| \leq L(x) \|\theta_{1} - \theta_{2}\|.$$
%Furthermore, assume that $P m_{\theta}$ admits a second-order Taylor expansion in $\theta$ at a point of maximum $\theta_{0}$ with nonsingular symmetric second derivative matrix $V_{\theta_{0}}$. If $P_{n} m_{\hat{\theta}_{M}} \geq \sup_{\theta} P_{n} m_{\theta} - o_{P}(n^{-1})$ and $\hat{\theta}_{M} \to^{P} \theta_{0}$ then $\sqrt{n}(\hat{\theta}_{M} - \theta_{0})$ is asymptotically normal with mean zero and covariance matrix $V_{\theta_{0}}^{-1} P \nabla{m}_{\theta_{0}} \nabla{m}_{\theta_{0}}^{\top} V_{\theta_{0}}^{-1}$.
%\end{theorem}
Interestingly, when applying this theorem in the context of maximum likelihood, differentiability in quadratic mean provides the existence of this second-order Taylor expansion. It is remarkable that a first-order assumption can be used to provide second-order properties; we refer the reader to \cite{Pollard} for further discussion on this point.

MLEs and M-estimators are often computed with iterative methods, so it is important to quantify the performance of estimators derived as iterates of an optimization method. One approach for doing so is one-step estimation, which was first introduced by Le Cam \cite{LeCamOG}. In its original form, one-step estimation consisted of applying one Newton-Raphson iteration, on a maximum log-likelihood objective, to an initial estimator within some range of the parameter to be estimated. Le Cam showed that, subject to certain conditions, this one-step estimator is asympotically equivalent to the MLE. This is remarkable because it shows that the result of a single Newton-Raphson iteration posesses the same statistical optimality as the limit point of Newton-Raphson iterates, the MLE. One-step estimation has since been extended to M-estimation \cite{bickel1975one}, sparse estimation \cite{zou2008, taddy2017one}, quasi-likelihood estimation \cite{fan1999one}, and distributed computation \cite{huang2019distributed}.

The one-step estimation result most relevant to our work is the following theorem, which summarizes the asymptotic performance of a certain one-step estimator.
\begin{theorem}[\cite{vanderVaart}, Theorem 5.21 and Theorem 5.45] \label{OSE_OG}
  Let $\hat{\theta}^{*}$ be an M-estimator with a continuously differentiable criterion $m_{\theta}$, where $E [\|\nabla m_{\theta^*}\|^2] <\infty$, $E[ \nabla m_{\theta^{*}}] = 0$, $E[m_{\theta}]$ is twice differentiable at $\theta^{*}$ with nonsingular Hessian $V_{\theta^{*}}$, and $\frac{1}{n} \sum_{i=1}^{n} \nabla m_{\hat{\theta}_{M}}(X_{i}) = o_{P}(n^{-1/2})$. Assume that $\hat{\theta}_{M}$ is consistent for $\theta^{*}$ and that there is a measurable function $L$ with $E[L^2] < \infty$ such that for every $\theta_{1}$ and $\theta_{2}$ in a neighborhood of $\theta^{*}$
  $$\left\| \nabla m_{\theta_{1}}(x) - \nabla m_{\theta_{2}}(x)\right\| \leq L(x) \left\|\theta_{1} - \theta_{2}\right\|.$$
  Let $\hat{\theta}_{init}$ be a $\sqrt{n}$-consistent estimator of $\theta^{*}$. If $C_{n}$ is a sequence of random matrices such that $C_{n} \to^{P} V_{\theta^{*}}$, then $\hat{\theta}_{M}$ and $\hat{\theta}_{init} - C_{n}^{-1} \left(\frac{1}{n} \sum_{i=1}^{n} \nabla m_{\hat{\theta}_{init}}(X_{i})\right)$ are asymptotically equivalent.
\end{theorem}

This result on one-step estimation gives that, subject to some mild conditions, one Newton step performed on an initial estimator within $n^{-1/2}$ of $\hat{\theta}_{M}$ gives the same large sample performance as the M-estimator itself. 

Though one-step estimation is defined as a single step of an iterative method, it can be interpreted as the \emph{last} step of an iterative method. As long as the penultimate iterate satisfies the properties required of $\hat{\theta}_{init}$, the last iteration generates an estimator with one-step properties. Moreover, further iterations do not increase the performance of the estimator from the asymptotic perspective. Theorem \ref{OSE_OG} demonstrates that the one-step estimator is asymptotically equivalent to $\hat{\theta}_{M}$; hence it is $\sqrt{n}$-consistent. The one-step estimator then can be used as $\hat{\theta}_{init}$ in another iteration, because it satisfies the required properties for the starting point, but the resulting asymptotic distribution remains the same. We conclude that numerical tolerance for Newton's method should be $O(n^{-1/2})$, where $n$ is the sample size of the problem, in order to respect the statistical origin of the problem and guarantee the asymptotic properties in Theorem \ref{OSE_OG}.

\subsection{Proximal Methods.}\label{proxmethods}

In this section we review proximal descent and proximal gradient descent, in addition to their extensions scaled proximal descent and proximal Newton descent. Each of these algorithms are iterative methods for finding local minima, and in each step these methods use a \emph{proximal operator}, a numerical primitive which involves solving a related optimization problem.

For a function $f: \R^{d} \to \R \cup \{\infty\}$ and positive parameter $\lambda > 0$, the proximal operator $\prox_{f}^{\lambda}(x)$ is defined as
$$\proxi{f}{\lambda}(x) = \argmin_{w \in \R^d}\left\{ f(w) + \frac{1}{2\lambda} \|w-x\|^2\right\}.$$
The proximal operator is a generalization of a projection. Indeed, if $S \subset \R^d$ and $I_{S}$ is an \emph{indicator function} on the set $S$
$$I_{S}(x) = \left\{\begin{array}{cc} 0 & \text{ if } x \in S\\
\infty & \text{ otherwise} \end{array}\right.$$
then the proximal operator of $I_{S}$ for any $\lambda >0$ is $\argmin_{w \in S} \|w-x\|^2$, projection onto the set $S$. Alternatively, one can view the proximal operator as a kind of implicit gradient descent. If $f$ is smooth and convex then the optimality conditions give
$$x - \lambda \nabla f\left(\proxi{f}{\lambda}(x)\right) = \proxi{f}{\lambda}(x).$$
Thus applying the proximal operator is equivalent to taking a gradient step, where the gradient is evaluated at the output of the proximal operator instead of the initial point $x$, and the step length is $\lambda$.

If the objective function $f$ is proper, so that it is finite for at least one point and never takes the value $-\infty$, then $x^{*}$ minimizing $f$ implies $\prox_{f}^{\lambda}(x^{*}) = x^{*}$. Since the proximal operator is also \emph{firmly nonexpansive} \cite[Prop.12.28]{BauCom11Cvx} when $f$ is convex, so that $\|\prox_{f}^{\lambda}(x) - \prox_{f}^{\lambda}(y)\| \leq \|x-y\|$, a natural approach to minimizing a convex $f$ is to iterate the proximal operator \cite[Prop.12.29]{BauCom11Cvx}. This algorithm is called \emph{proximal descent}, and it converges for any sequence $\lambda_k$ of the $\lambda$ parameters such that $\lambda >0$ and $\sum_{k} \lambda_k \to \infty$. Proximal descent is notable for its simplicity, but it finds limited application in practice because its convergence rate matches gradient descent \cite{BeckBook}. Evaluation of the proximal operator is rarely as simple as gradient evaluation, so besides a few notable exceptions \cite{mattingley2012cvxgen, golub1966note} proximal descent is primarily of theoretical interest. Additional details of the proximal operator and the proximal descent method can be found in \cite[\S12.4, \S28.5]{BauCom11Cvx}.

Other methods based on proximal descent are extremely useful in practice. Consider the composite model
$$\min_{x \in \R^d} \; g(x) + h(x)$$
where $g: \R^d \to \R$ is closed and continuously differentiable and $h: \R^d \to \R \cup \{\infty\}$ is convex and potentially nonsmooth. If a point $x^{*}$ minimizes $g(x) + h(x)$ then it can easily be shown that\footnote{using, for example, Theorem 10.1 and Proposition 8.12 in \cite{VarAn}.}
%Using Theorem 10.1, Exercise 8.8, and Proposition 8.12 in Variational Analysis

%this discussion follows Beck and Teboulle. I took it from Boyd's Proximal Algorithm Treatise
$$0 \in \nabla g(x^{*}) + \partial h(x^{*})$$
where $\partial$ denotes the subdifferential of a convex function. A few algebraic manipulations give that
\begin{equation*}\label{eq:fixedpt}
  x^{*} - \lambda \nabla g(x^{*}) \in x^{*} + \lambda \partial h(x^{*}).
\end{equation*}
Hence we conclude that $x^{*} \in \prox_{h}^{\lambda}\left( x^{*} - \lambda \nabla g(x^{*})\right)$. If in addition $g$ is strongly convex function with $L$-Lipschitz continuous gradient the operator $\prox_{h}^{\lambda}\left( x - \lambda \nabla g(x)\right)$ is firmly nonexpansive for values of $\lambda$ not larger than $\frac{1}{L}$ \cite{PaBo14proximal}. Thus, we iterate  $\prox_{h}^{\lambda}\left( x - \lambda \nabla g(x)\right)$ in order to find a minimizer of $g(x) + h(x)$. Because each step of this algorithm composes a proximal step in $h$ with a gradient descent step in $g$, this method is called \emph{proximal gradient descent}. Proximal gradient descent finds widespread application in many problems in statistics and machine learning, where the strongly convex function $g$ is most often a data fidelity term and $h$ is a penalty or constraint which encourages certain solution structure \cite{polson2015proximal}. In the setting of Bayesian point estimation, the function $h$ can be used to incorporate a prior using Maximum-a-Posteriori (MAP) estimation \cite{BassettDeride}.

As alternatives to proximal descent and proximal gradient descent, we consider their scaled analogues. For a function $f$ and positive definite matrix $C$ we define the \emph{scaled proximal operator} as 
$$\proxi{f}{C}(x) := \argmin_{w \in \R^d}\left\{f(w) + \frac{1}{2} \|w - x\|^{2}_{C}\right\},$$
where $\|\cdot\|_C$ is the scaled norm induced by $C$, i.e., $\|x\|^2_C=\langle x,\,Cx\rangle$ for every $x$ in $\R^d$. By replacing the proximal operator with the scaled proximal operator in proximal descent, we have scaled proximal descent. Similarly, replacing the step length $\lambda$ in proximal gradient descent with a scaling matrix $C^{-1}$ yields scaled proximal gradient descent. The only notable difference is that the firm nonexpansivity in both cases is now with respect to the $\|\cdot\|_{C}$ norm. Various conventions to generate $C$ have been previously investigated, including scalar multiples of the identity \cite{BeTe09fista,Mil16numerical}, diagonal matrices \cite{Tse09coordinate}, quasi-Newton approximations of $\nabla^2 f$ \cite{sr1, lbfgs,Kan21globalized}, and taking $C = \nabla^{2} f$ \cite{leesunsaunders}. Better approximations of the curvature of $f$ generally make the scaled proximal operator harder to evaluate, but also result in scaled proximal descent requiring fewer iterations to reach a fixed solution accuracy \cite{FriedlanderGoh}.

%\textbf{This block should probably be purged in a final version}
%
%Indeed, let $x$ and $y$ be points in $\R^d$. Let $w = \prox_{f, C}(x)$ and $v = \prox_{f, C}(y)$. Then by the Cauchy-Schwarz inequality and the monotonicity of a convex subdifferential
%\begin{align*}
%\|x-y\|_{C^{-1}} \|w-v\|_{C^{-1}} & \geq (x-y)^{\top} C^{-1} (w-v) \\
%& = \left(w  + C \partial h(w) - (v + C \partial h(v))\right)^{\top} C^{-1} (w-v)\\
%& = \|w -v \|^{2}_{C^{-1}} + \left(\partial h(w) - \partial h(v)\right)^{\top} (w-v)\\
%& \geq \|w - v\|^{2}_{C^{-1}}
%\end{align*}
%Dividing by $\|w-v\|_{C^{-1}}$ on each side gives firm nonexpansivity of the scaled proximal operator in the $C^{-1}$ norm.
%
%Similarly, a descent step with scaling matrix $C$ satisfies
%\begin{align*}
%\|x -y\|_{C^{-1}}^2 &= (x -y)^{\top} C^{-1} (x -y) \\
%& \leq (x-y)^{\top} C^{-1} (x - C \nabla f(x) - (y - C \nabla f(y)))\\
%& \leq \| x- y\|_{C^{-1}} \left\|x - C\nabla f(x) - (y - C\nabla f(y))\right\|_{C^{-1}}.
%\end{align*}
%Dividing by $\| x- y\|_{C^{-1}}$ on both sides gives firm nonexpansivity of a scaled descent step in the $C^{-1}$ norm. Combining firm nonexpansivity of this descent step and the proximal step gives firm nonexpansivity of a proximal Newton iteration in the $C^{-1}$ norm.
%
%\textbf{End Block}

\section{The Composite Model.}\label{Composite} 

A composite model
%\footnote{Calling $g+h$ where $g$ is strongly convex and $g$ is convex the \emph{composite model} is found in \cite{BeckBook}[Chapter 10].} 
refers to an optimization problem $\min f$, where $f = g+h$. It is common practice to assume that $g$ is convex with Lipschitz continuous gradient whereas the function $h$ is convex and potentially nonsmooth. Where possible, we will also relax the convexity assumption on $g$.

In this section, we study the asymptotic equivalence of an estimator defined as the minimizer in a composite model and the one-step estimator attained through scaled proximal gradient descent. Let $g_{n}: \Omega \times \Theta \to \mathbb{R}$ and  $h_{n}: \Omega \times \Theta \to \mathbb{R} \cup \{\infty\}$ be sequences of random functions on $\Theta \subseteq \mathbb{R}^{d}$, defined on a probability space $(\Omega, \mathcal{A}, P)$, and let $C_{n}$ be a sequence of $d\times d$ random matrices defined on the same probability space. Note that this setting is more general than the parametric estimation from i.i.d. observations introduced in subsection \ref{subsec:one-step}. We will return to the i.i.d. setting in Proposition \ref{prop:assump1}.

%In what follows, we study the asymptotic equivalence of an estimator defined as the minimizer of the composite model and the one-step estimator attained through scaled proximal gradient descent. Let $g_{n}: \Omega \times \Theta \to \mathbb{R}$ and  $h_{n}: \Omega \times \Theta \to \mathbb{R} \cup \{\infty\}$ be sequences of random criteria functions on $\Theta$ defined on a common probability space $(\Omega, \mathcal{F}, P)$. 

We suppress the dependence on $\Omega$ in our notation, so that $g_{n}(\theta)$ and $h_{n}(\theta)$ denote a sequence of random variables indexed by $\theta \in \Theta$. Similarly, gradients of $g_{n}$ and subgradients of $h_{n}$ are with respect to the $\Theta$ argument. Define an estimator for the composite model $g_{n} + h_{n}$ as
\begin{equation}\label{eq:compmodel}
\hat{\theta} = \argmin_{\theta\in\Theta} \; g_{n} + h_{n},
\end{equation}
which is assumed to exist $P$-almost surely. Note that $\hat{\theta}$ is a function of $n$, though to ease notation we have omitted this dependence.

For an initial estimator $\hat{\theta}_{init}$ of $\hat{\theta}$, define the associated one-step estimator (OSE), denoted $\hat{\theta}_{ose}$, as one iteration of scaled proximal gradient descent from the initial estimator $\hat{\theta}_{init}$,
\begin{equation}\label{def:cmose}
\hat{\theta}_{ose} = \proxi{h_{n}}{C_{n}}(\hat{\theta}_{init} - C_{n}^{-1} \nabla g_{n}(\hat{\theta}_{init})),
\end{equation}
which is also assumed to exist $P$-almost surely. The initial estimator $\hat{\theta}_{init}$ is assumed to be a $\sqrt{n}$-consistent estimator of $\hat{\theta}$, so that $\sqrt{n}(\hat{\theta}_{init} - \hat{\theta})$ is bounded in probability.

Our first result provides conditions under which the estimators in equations \eqref{eq:compmodel} and \eqref{def:cmose} are asymptotically equivalent, so that the one-step estimator has the same asymptotic performance as the minimizer of the composite model.

%{\color{red} It's important that we do not need $g_{n}$ convex, otherwise this invalidates our Cauchy likelihood example. However, when $g_{n}$ is not convex we need to be extra clear that a minimizer is a fixed point of proximal gradient descent. Exercise 8.8 in Variational Analysis may be useful in that effort.}

\begin{theorem} \label{thm:proxgrad}
  Consider the composite estimation model in equation~(\ref{eq:compmodel}), where it holds with probability 1 that $g_{n}$ is continuously differentiable and $h_{n}$ is convex and proper.

Assume the following.
\begin{enumerate}[(i)]
\item For each $M>0$, there is a positive definite matrix $H$ such that
  $$\sup_{\sqrt{n}\norm{\theta-\hat{\theta}} < M} \sqrt{n}\left(\nabla g_{n}(\theta) - \nabla g_{n}(\hat{\theta})\right) - \sqrt{n} H (\theta - \hat{\theta}) \to^{P} 0.$$
\item $C_n\succ 0$ and $C_{n}^{-1} H \to^{P} I$, where $I$ denotes the $d \times d$ identity matrix.
%\item With high probability there exists some constant $K$ such that $\prox_{h_{n}}^{C_{n}}$ is lipschitz continuous with constant $K$.
%\item $\hat{\theta} - \prox_{h_{n}}^{C_{n}}\left(\hat{\theta} - C_{n}^{-1} \nabla g_{n}(\hat{\theta})\right) = o_{p}(n^{-1/2})$. That is, $\hat{\theta}$ is approximately a fixed point of scaled proximal gradient descent.
\end{enumerate}
  Then for any $\sqrt{n}$-consistent estimator $\hat{\theta}_{init}$ of $\hat{\theta}$, $\hat{\theta}$ and the one-step estimator $\hat{\theta}_{ose}$ of equation~(\ref{def:cmose}) are asymptotically equivalent.
\end{theorem}

\proof{}
  We first show that for some constant $K$, the scaled proximal operator $\prox_{h_{n}}^{C_{n}}$ is $K$-Lipschitz continuous with high probability. Indeed, for any $\theta$ the proximal operator satisfies
  $$C_{n} \left(\theta - \proxi{h_{n}}{C_{n}}(\theta)\right) \in \partial h \left(\proxi{h_{n}}{C_{n}}(\theta)\right).$$
  So by monotonicity of the subdifferential of $h_{n}$, for any $\theta_{1}$, $\theta_{2}$ we have
  $$\left(\proxi{h_{n}}{C_{n}}(\theta_{1}) - \proxi{h_{n}}{C_{n}}(\theta_{2})\right)^{T} C_{n} \left(\theta_{1} - \proxi{h_{n}}{C_{n}}(\theta_{1}) - \left(\theta_{2} - \proxi{h_{n}}{C_{n}}(\theta_{2}) \right) \right) \geq 0.$$
  Applying the Cauchy-Schwarz inequality gives
  $$\left\|\theta_{1} - \theta_{2}\right\|_{C_{n}} \geq \left\|\proxi{h_{n}}{C_{n}}(\theta_{1}) - \proxi{h_{n}}{C_{n}}(\theta_{2}) \right\|_{C_{n}}$$
  Therefore $\prox_{h_{n}}^{C_{n}}$ is firmly nonexpansive in the $C_{n}$ norm. Denote by $\lambda_{max}(C_{n})$ and $\lambda_{min}(C_{n})$ the maximum and minimum eigenvalues of $C_{n}$. We have
  $$\sqrt{\frac{\lambda_{max}(C_{n})}{\lambda_{min}(C_{n})}} \left\| \theta_{1} - \theta_{2}\right\| \geq \left\|\proxi{h_{n}}{C_{n}}(\theta_{1}) - \proxi{h_{n}}{C_{n}}(\theta_{2}) \right\|.$$
  Because $C_{n}^{-1} H \to^{P} I$ and $H \succ 0$, $\sqrt{\frac{\lambda_{max}(C_{n})}{\lambda_{min}(C_{n})}}$ converges in probability to some constant. It follows that there exists a $K$, for any $\epsilon > 0$, and for $n$ large enough such that  $\sqrt{\lambda_{max}(C_{n})/\lambda_{min}(C_{n})} < K$ with probability $1-\epsilon$. We conclude that with high probability $\prox_{h_{n}}^{C_{n}}$ is $K$-Lipschitz continuous.  

We next want to show $\sqrt{n} \norm{\hat{\theta}_{ose} - \hat{\theta}} = o_{P}(1)$. From the definition of $\hat{\theta}_{ose}$ we have
  \begin{align*}
    \sqrt{n} \left\| \hat{\theta}_{ose} - \hat{\theta} \right\|
    \leq &\sqrt{n} \left\| \proxi{h_{n}}{C_{n}}\left(\hat{\theta}_{init} - C_{n}^{-1} \nabla g_{n}(\hat{\theta}_{init}) \right) - \proxi{h_{n}}{C_{n}}\left(\hat{\theta} - C_{n}^{-1} \nabla g_{n}(\hat{\theta}) \right)\right\| \\& +  \sqrt{n} \left\|\proxi{h_{n}}{C_{n}}\left(\hat{\theta} - C_{n}^{-1} \nabla g_{n}(\hat{\theta}) \right) - \hat{\theta}\right\|
  \end{align*}
  Recall from section \ref{proxmethods} that $\hat{\theta}$ is a fixed point of scaled proximal gradient descent. Therefore the final term in the sum is $o_{P}(1)$, and we can simplify as follows.
  \begin{align*}
    = & \sqrt{n} \left\| \proxi{h_{n}}{C_{n}}\left(\hat{\theta}_{init} - C_{n}^{-1} \nabla g_{n}(\hat{\theta}_{init}) \right) - \proxi{h_{n}}{C_{n}}\left(\hat{\theta} - C_{n}^{-1} \nabla g_{n}(\hat{\theta}) \right)\right\| + o_{P}(1). 
  \end{align*}
  From the $\sqrt{n}$-consistency of $\hat{\theta}_{init}$ and Lipschitz continuity result above we can restrict our considerations to the event where $\sqrt{n} \| \hat{\theta}_{init} - \hat{\theta} \| < M $ for some $M$ and $\prox_{h_{n}}^{C_{n}}$ is Lipschitz continuous with constant $K$, since its complement can be made to have arbitrarily small probability. On this event, assumption $(i)$ holds and the above display becomes the following.
  \begin{align}
    \leq & K \sqrt{n} \left\| \hat{\theta}_{init} - C_{n}^{-1} \nabla g_{n}(\hat{\theta}_{init}) - \left(\hat{\theta} - C_{n}^{-1} \nabla g_{n}(\hat{\theta}) \right)\right\| + o_{P}(1) \nonumber\\
    = & K \sqrt{n} \left\| \left(I - C_{n}^{-1} H \right) \left(\hat{\theta}_{init} - \hat{\theta}\right) + C_{n}^{-1} H (\hat{\theta}_{init} - \hat{\theta}) - C_{n}^{-1} \left(\nabla g_{n}(\hat{\theta}_{init}) - \nabla g_{n}(\hat{\theta}) \right) \right\| + o_{P}(1)\nonumber\\
%    \leq & K \sqrt{n} \left\|\left(I - C_{n}^{-1} H\right) \left(\hat{\theta}_{init} - \hat{\theta}\right) \right\| + K \sqrt{n} \left\| C_{n}^{-1} \left( \nabla g_{n}(\hat{\theta}_{init}) - \nabla g_{n}(\hat{\theta}) - H(\hat{\theta}_{init} - \hat{\theta})\right) \right\| + o_{P}(1)\nonumber\\
%    \leq & K \sqrt{n} \left\|I - C_{n}^{-1} H\right\| \left\|\hat{\theta}_{init} - \hat{\theta} \right\| + K \sqrt{n} \left\| C_{n}^{-1} \right\| \left\| \nabla g_{n}(\hat{\theta}_{init}) - \nabla g_{n}(\hat{\theta}) - H(\hat{\theta}_{init} - \hat{\theta}) \right\| + o_{P}(1)\nonumber\\
    \leq & K \left\|I - C_{n}^{-1} H\right\| \cdot \sqrt{n} \left\|\hat{\theta}_{init} - \hat{\theta} \right\| + K  \left\| C_{n}^{-1} \right\| \cdot \sqrt{n} \left\| \nabla g_{n}(\hat{\theta}_{init}) - \nabla g_{n}(\hat{\theta}) - H(\hat{\theta}_{init} - \hat{\theta}) \right\| + o_{P}(1) \label{eq:almostdone}
  \end{align}
    From assumption $(ii)$, $\left\|I - C_{n}^{-1} H \right\| \to^{P} 0$. Because $\hat{\theta}_{init}$ is $\sqrt{n}$-consistent for $\hat{\theta}$, $\sqrt{n}\left\| \hat{\theta}_{init} - \hat{\theta}\right\| = O_{P}(1)$. Further, $K \|C_{n}^{-1}\| = O_{P}(1)$ because of assumption $(ii)$. Lastly, assumption $(i)$ gives that 
    $$\sqrt{n} \left\| \nabla g_{n}(\hat{\theta}_{init}) - \nabla g_{n}(\hat{\theta}) - H(\hat{\theta}_{init} - \hat{\theta}) \right\| = o_{P}(1).$$
    Therefore the expression in \eqref{eq:almostdone} is $o_{P}(1)$ by Slutsky's Theorem. This concludes the proof.
\endproof

We comment briefly on the assumptions in Theorem~\ref{thm:proxgrad}. The convexity of $h$ was only used to establish the Lipschitz continuity of the scaled proximal operator. By restricting to the event $\{\sqrt{n}\|\hat{\theta}_{init} - \hat{\theta}\|\}$, it suffices to allow $h$ to belong to a wider class of functions, such as lower-semicontinuous and prox-regular, where it can be established that the scaled proximal operator is locally Lipschitz continuous \cite[Theorem 1]{HareSagastizabal}.

Assumption $(i)$ is analogous to a similar condition in \cite{vanderVaart}, where it is used to demonstrate asymptotic optimality of maximum likelihood estimators and one-step estimators. Despite its resemblance to a condition guaranteeing $\nabla g_{n}$ is differentiable with Jacobian $H$, the probabilistic nature of the limit allows it to be satisfied even when $\nabla g_{n}$ is not differentiable. In this case, applying the well known {\it discretization trick} gives that differentiability in quadratic mean implies assumption $(i)$ when $\hat{\theta}$ converges in probability to some nonrandom $\theta_{0}$ (see \cite{cam} for details). %\cite{vanderVaart}.

Assumption $(ii)$ is slightly stronger than necessary; we only require that 
$$\sqrt{n} \left(I - C_{n}^{-1} H\right) \left(\hat{\theta}_{init} - \hat{\theta}\right) \to^{P} 0.$$
This assumption has an interesting relationship with the Dennis-Mor\'e criterion, a condition for the approximation of Hessians in deterministic variable metric methods \cite{dennismore}. The Dennis-Mor\'e criterion states that the approximation $C_{n}$ of a Hessian $H$ satisfies
$$\frac{\left\|\left(C_{n} - H\right) \left(\hat{\theta}_{k+1} - \hat{\theta}_{k}\right) \right\|}{\left\|\hat{\theta}_{k+1} - \hat{\theta}_{k}\right\|} \to 0,$$
where $\hat{\theta}_{k+1}$ and $\hat{\theta}_{k}$ are iterates of a variable metric method. Assumption $(ii)$ is implied by a probabilistic version of the Dennis-Mor\'e criterion, where iterates are replaced by the $\hat{\theta}_{init}$ and $\hat{\theta}$. Indeed, if 
$$\frac{\left\|\left(C_{n} - H\right) \left(\hat{\theta}_{init} - \hat{\theta}\right) \right\|}{\left\|\hat{\theta}_{init} - \hat{\theta}\right\|} = o_{p}(1)$$
then using the $\sqrt{n}$-consistency of $\hat{\theta}_{init}$
$$\sqrt{n} \left\|\left(C_{n} - H\right) \left(\hat{\theta}_{init} - \hat{\theta}\right) \right\|
= O\left(\frac{\left\|\left(C_{n} - H\right) \left(\hat{\theta}_{init} - \hat{\theta}\right) \right\|}{\left\|\hat{\theta}_{init} - \hat{\theta}\right\|}\right) = o_{p}(1),$$
yielding assumption $(ii)$ whenever $\left\|C_{n}\right\|$ is bounded in probability.

We also note that $\hat{\theta}$ need not be a minimizer of $g_{n} + h_{n}$ for the conclusion of Theorem \ref{thm:proxgrad} to hold. Indeed, the only required property of $\hat{\theta}$ is that it is approximately a fixed point of scaled proximal gradient descent. That is,
\begin{equation}
\label{eq:approxmin}
\sqrt{n} \left\|\proxi{h_{n}}{C_{n}}\left(\hat{\theta} - C_{n}^{-1} \nabla g_{n}(\hat{\theta}) \right) - \hat{\theta}\right\| = o_{P}(1).
\end{equation}
This flexibility allows Theorem \ref{thm:proxgrad} to generalize to stationary points, because $g_{n}$ need not be convex. We formalize this result in the following corollary.

\begin{corollary} \label{cor:approxmin}
  In the setting of Theorem \ref{thm:proxgrad}, if $\hat{\theta}$ satisfies \eqref{eq:approxmin} as an approximate fixed point of proximal gradient descent, then the conclusion of Theorem~\ref{thm:proxgrad} holds even if $\hat{\theta}$ is not necessarily a minimizer of $g_{n} + h_{n}$.
%  is , but instead satisfies \eqref{eq:approxmin} 
\end{corollary}

We conclude this subsection with a result which permits the application of Theorem \ref{thm:proxgrad} to regularized maximum likelihood in the setting of i.i.d. observations from subsection \ref{subsec:one-step}. The following result establishes the asymptotic equivalence of the regularized maximum likelihood estimator and the one-step estimator derived from scaled proximal gradient descent. Its assumptions are classical and verifiable for a large class of parametric models \cite[\S5.6]{vanderVaart}.

\begin{proposition} \label{prop:assump1}
  Assume that there is some $\theta_{0} \in \Theta$ such that the estimator $\hat{\theta}$ in \eqref{eq:compmodel} is $\sqrt{n}$-consistent for $\theta_{0}$. Assume also that $X_{i} \sim^{i.i.d.} p_{\theta^{*}}$ for some $\theta^{*} \in \Theta$. Let $g_{n}(\theta) = - \frac{1}{n} \sum_{i=1}^{n} \log p_{\theta}(X_{i})$ and $h_n$ be proper and convex $P$-almost surely.  Assume that
  \begin{enumerate}[(i)]
    \item $\log p_{\theta}$ is three-times differentiable for every $x$.
    \item The matrix $E[\nabla^{2} \log p_{\theta_{0}}]$ exists and is negative definite (so that the Fisher Information $I_{\theta_{0}}$ exists and is nondegenerate).
    \item The third order partial derivatives of $\log p_{\theta}$ are dominated by a fixed integrable function for every $\theta$ in a neighborhood of $\theta_{0}$.
    \item $C_n\succ 0$ and $-C_n^{-1}E[\nabla^{2} \log p_{\theta_{0}}] \to^{P} I$. \label{assump1:4}
\end{enumerate}
%\item With high probability there exists some constant $K$ such that $\prox_{h_{n}}^{C_{n}}$ is lipschitz continuous with constant $K$.
%\item $\hat{\theta} - \prox_{h_{n}}^{C_{n}}\left(\hat{\theta} - C_{n}^{-1} \nabla g_{n}(\hat{\theta})\right) = o_{p}(n^{-1/2})$. That is, $\hat{\theta}$ is approximately a fixed point of scaled proximal gradient descent.
%\end{enumerate}
  Then for any $\sqrt{n}$-consistent estimator $\hat{\theta}_{init}$ of $\hat{\theta}$, $\hat{\theta}$ and the one-step estimator $\hat{\theta}_{ose}$ of equation~(\ref{def:cmose}) are asymptotically equivalent.
%  Then Theorem \ref{thm:proxgrad} holds with only assumption $(ii)$.
\end{proposition}

\proof{}
  The result follows from applying Theorem~\ref{thm:proxgrad} under these assumptions.  Thus, we need to show that for each $M>0$ there is positive definite $H$ such that
  $$\sup_{\sqrt{n}\norm{\theta-\hat{\theta}} < M} \sqrt{n}\left(\nabla g_{n}(\theta) - \nabla g_{n}(\hat{\theta})\right) - \sqrt{n} H (\theta - \hat{\theta}) \to^{P} 0.$$
Take $H = -E[\nabla^{2} \log p_{\theta_{0}}]$, so that assumption $(iv)$ is equivalent to Theorem~\ref{thm:proxgrad}'s assumption $(ii)$. Fix $M > 0$. By $\sqrt{n}$-consistency of $\hat{\theta}$, $\sqrt{n}\norm{\hat{\theta} - \theta_{0}} < K$ for some $K$ since its complement can be made to have arbitrarily small probability. We thus have $\sqrt{n}\norm{\theta - \theta_{0}} < M + K$ when $\sqrt{n} \left\|\theta - \hat{\theta}\right\| < M$. We can expand with the triangle inequality to give
  \begin{align}
    & \left\|\sup_{\sqrt{n}\norm{\theta-\hat{\theta}} < M} \sqrt{n}\left(\nabla g_{n}(\theta) - \nabla g_{n}(\hat{\theta})\right) - \sqrt{n} H (\theta - \hat{\theta})\right\| \nonumber \\
    & \leq \sup_{\sqrt{n}  \norm{\theta - \hat{\theta}} < M} \bigg\{ \norm{\sqrt{n} \left( \nabla g_{n}(\theta) - \nabla g_{n}(\theta_{0}) \right) - \sqrt{n} \nabla^{2} g_{n}(\theta_{0})(\theta - \theta_{0})}  \nonumber\\
    &  + \norm{\sqrt{n} \left( \nabla g_{n}(\theta_{0}) - \nabla g_{n}(\hat{\theta}) \right) - \sqrt{n} \nabla^{2} g_{n}(\theta_{0})(\theta_{0} - \hat{\theta})} \nonumber\\
    & + \norm{\sqrt{n}\left(\nabla^{2} g_{n}(\theta_{0}) - H \right)(\theta - \theta_{0})} + \norm{\sqrt{n} \left(\nabla^{2} g_{n}(\theta_{0}) - H \right)(\theta_{0} - \hat{\theta})} \bigg\}\nonumber\\
    & \leq \sup_{\sqrt{n}  \norm{\theta - \theta_{0}} < M + K} \norm{\sqrt{n} \left( \nabla g_{n}(\theta) - \nabla g_{n}(\theta_{0}) \right) - \sqrt{n} \nabla^{2} g_{n}(\theta_{0})(\theta - \theta_{0})} \nonumber\\
    & + \norm{\sqrt{n} \left( \nabla g_{n}(\theta_{0}) - \nabla g_{n}(\hat{\theta}) \right) - \sqrt{n} \nabla^{2} g_{n}(\theta_{0})(\theta_{0} - \hat{\theta})} \label{eq:bigsum}\\
    & + \sup_{\sqrt{n}\norm{\theta - \theta_{0}} < M+K} \norm{\sqrt{n}\left(\nabla^{2} g_{n}(\theta_{0}) - H \right)(\theta - \theta_{0})} \nonumber\\
    & + \norm{\sqrt{n} \left(\nabla^{2} g_{n}(\theta_{0}) - H \right)(\theta_{0} - \hat{\theta})} \nonumber
  \end{align}
  We have $\nabla^{2} g_{n}(\theta_{0}) = -\frac{1}{n}\sum_{i=1}^{n} \nabla^{2} \log p_{\theta_{0}}(X_{i})$, which converges almost surely to $H$ by the law of large numbers. Since both $\sqrt{n}(\theta - \theta_{0})$ and $\sqrt{n}(\hat{\theta} - \theta_{0})$ are bounded in probability, the third and fourth expressions in \eqref{eq:bigsum} are $o_{P}(1)$. We next focus on the first term in the \eqref{eq:bigsum}. We have
  $$\sqrt{n} \left(\nabla g_{n}(\theta) - \nabla g_{n}(\theta_{0})\right) - \sqrt{n} \nabla^{2} g_{n}(\theta_{0})(\theta - \theta_{0})$$
  $$= -\frac{1}{\sqrt{n}} \sum_{i=1}^{n} \nabla \log p_{\theta}(X_{i}) - \nabla \log p_{\theta_{0}}(X_{i}) - \nabla^{2} \log p_{\theta_{0}}(X_{i}) (\theta - \theta_{0})$$
  By Taylor's Theorem and assumption $(iii)$, 
  $$\left\|\nabla \log p_{\theta}(X_{i}) - \nabla \log p_{\theta_{0}}(X_{i}) - \nabla^{2} \log p_{\theta_{0}}(X_{i}) (\theta - \theta_{0}) \right\| \leq C(X_{i}) \left\|\theta - \theta_{0}\right\|^{2}$$
  for some integrable function $C$. Therefore 
  $$\left\| \sqrt{n} \left(\nabla g_{n}(\theta) - \nabla g_{n}(\theta_{0})\right) - \sqrt{n} \nabla^{2} g_{n}(\theta_{0})(\theta - \theta_{0}) \right\|$$
  $$\leq n^{-1/2} \left\| \theta - \theta_{0} \right\|^{2} \sum_{i=1}^{n} C(X_{i})$$
  $$= n \left\| \theta - \theta_{0} \right\|^{2} n^{-3/2} \sum_{i=1}^{n}  C(X_{i})$$
  Since $n \left\| \theta - \theta_{0}\right\|^{2}$ is bounded in the supremum and $\frac{1}{n} \sum_{i=1}^{n} C(X_{i})$ is bounded in probability by the law of large numbers, we conclude that
  $$\sup_{\sqrt{n}  \norm{\theta - \theta_{0}} < M + K} \norm{\sqrt{n} \left( \nabla g_{n}(\theta) - \nabla g_{n}(\theta_{0}) \right) - \sqrt{n} \nabla^{2} g_{n}(\theta_{0})(\theta - \theta_{0})} \to^{P} 0.$$
  The term
  $$\norm{\sqrt{n} \left( \nabla g_{n}(\theta_{0}) - \nabla g_{n}(\hat{\theta}) \right) - \sqrt{n} \nabla^{2} g_{n}(\theta_{0})(\theta_{0} - \hat{\theta})}$$
   is $o_{P}(1)$ by similar reasoning, where we use the $\sqrt{n}$-consistency of $\hat{\theta}$ instead of the $\sqrt{n}  \norm{\theta - \theta_{0}} < M + K$ condition in the supremum. We conclude that
  $$\left\|\sup_{\sqrt{n}\norm{\theta-\hat{\theta}} < M} \sqrt{n}\left(\nabla g_{n}(\theta) - \nabla g_{n}(\hat{\theta})\right) - \sqrt{n} H (\theta - \hat{\theta})\right\| = o_{P}(1)$$
  which completes the proof.
\endproof

In Proposition \ref{prop:assump1}, we can take $\hat{\theta}_{init}$ to be any $\sqrt{n}$-consistent estimator of $\theta_{0}$ (which makes it $\sqrt{n}$-consistent for $\hat{\theta}$), such as a moment estimator. We emphasize that Proposition \ref{prop:assump1} permits $\theta_{0} \neq \theta^{*}$, which is of interest in structured inference problems where having estimates with certain properties (i.e. sparsity) may be more important than converging to the truth.

\subsection{Stopping Condition.}

The following result provides a recipe for finding a $\sqrt{n}$-consistent estimator $\hat{\theta}_{init}$. Given mild assumptions on the composite function $f_{n} = g_{n} + h_{n}$, $\hat{\theta}_{init}$ can be generated from a number of scaled proximal gradient steps and an easily verifiable stopping condition. Its proof is in the appendix.

\begin{proposition} \label{stop_cond}
  Let $f_{n} = g_{n} + h_{n}$ be a composite function and $\hat{\theta}$ a minimizer of $f_{n}$. Assume $h_{n}$ is proper and convex $P$-almost surely. Assume further that, with high probability, the function $g_n$ is strongly convex with parameter $m$ and $\nabla g_{n}$ is $M$-Lipschitz continuous. Let $C_{n}\succ 0$ be the matrix used in the scaled proximal gradient step. Assume that there is a constant $L$ such that $C_{n} \preceq L I$ with high probability. Then for any sequence $r_n$ such that $\|\hat{\theta}_{ose} - \hat{\theta}_{init} \| = O_{P}(r_{n})$, we have $\| \hat{\theta}_{init} -  \hat{\theta}\| = O_{P}(r_{n})$.
\end{proposition}

Proposition \ref{stop_cond} allows us to generate $\sqrt{n}$-consistent estimators from iterations of scaled proximal gradient descent. In practice, one can use an off-the-shelf implementation of scaled proximal gradient descent, such as the popular R package glmnet for penalized generalized linear models \cite{glmnet}, and take $\hat{\theta}_{init}$ to be any point for which the step length is less than $c/\sqrt{n}$ for some chosen constant $c$. We also note that Proposition \ref{stop_cond} applies when with high probability $\hat{\theta}$, $\hat{\theta}_{init}$, $\hat{\theta}_{ose}$ all lie in some set where the Lipschitz and strong convexity conditions hold.

\section{Proximal Descent.} \label{ProxDesc}

In this section we consider one-step estimators formed by minimizing an objective with scaled proximal descent. In contrast to the previous section, the objective function is written as a single function instead of the sum in the composite model. Consider the estimator
\begin{equation}\label{eq:estproxd}
\hat{\theta} = \argmin_{\theta} f_{n}(\theta),
\end{equation}
where as in section \ref{Composite}, $f_n: \Omega \times \Theta \to \mathbb{R} \cup \{\infty\}$ is a sequence of random functions on $\Theta \subseteq \mathbb{R}^{d}$, defined on a probability space $(\Omega, \mathcal{A}, P)$. For a sequence of random $d \times d$ scaling matrices $C_{n}$ and an initial estimator $\hat{\theta}_{init}$, we define a one-step estimator from proximal descent as
\begin{equation}\label{eq:oseproxd}
\hat{\theta}_{ose} = \proxi{f_{n}}{C_{n}}(\hat{\theta}_{init}).
\end{equation}
Compared to the scaled proximal gradient descent method in the previous section, this one-step estimator is a scaled proximal descent method (i.e. it omits the gradient step). We next provide a result which shows that $\hat{\theta}_{ose}$ in \eqref{eq:oseproxd} is asymptotically equivalent to $\hat{\theta}$ from \eqref{eq:estproxd}. This result is an analogue of Theorem \ref{thm:proxgrad} for scaled proximal descent.

\begin{theorem}\label{thm:proxdescent}
  Consider the estimator $\hat{\theta}$ in equation \eqref{eq:estproxd}, and $\hat{\theta}_{ose}$ in equation \eqref{eq:oseproxd}, where $f_{n}$ is such that both minimizers exist almost surely. Assume the following.
\begin{enumerate}[(i)]
\item $\lambda_{max}(C_{n}) = o_{P}(1)$ 
\item There is a constant $L$ such that for each constant $M$, the function $\prox_{f_{n}}^{C_{n}}$ is L-Lipschitz continuous on $\left\{\theta: \|\theta - \hat{\theta}\| \leq M/\sqrt{n}\right\}$ with high probability.
\item There is a constant $m$ such that for each constant $M$, the function $f_{n}$ is strongly convex with modulus $m$ on $\left\{\theta: \|\theta - \hat{\theta}\| \leq M/\sqrt{n}\right\}$ with high probability.
\end{enumerate}
  Then for any $\sqrt{n}$-consistent estimator $\hat{\theta}_{init}$ of $\hat{\theta}$, $\hat{\theta}$ and $\hat{\theta}_{ose}$ are asymptotically equivalent.
\end{theorem}

\proof{}
  As a global minimizer of $f_{n}$, $\hat{\theta}$ is a fixed point of $\prox_{f_{n}}^{C_{n}}$. Assume that $\sqrt{n}\left\|\hat{\theta} - \hat{\theta}_{init}\right\| < M$ for some constant $M$, since its complement has arbitrarily small probability by the $\sqrt{n}$-consistency of $\hat{\theta}_{init}$.
  
  We need to show that $\sqrt{n} \left\|\hat{\theta}_{ose} - \hat{\theta}\right\| = o_{P}(1)$. From the definition of $\hat{\theta}_{ose}$ we have the following.
  \begin{align}
    \sqrt{n} \left\|\hat{\theta}_{ose} - \hat{\theta}\right\| \leq & \sqrt{n}\left\|\hat{\theta}_{ose} - \proxi{f_{n}}{C_{n}}(\hat{\theta})\right\| + \sqrt{n}\left\|\proxi{f_{n}}{C_{n}}(\hat{\theta}) - \hat{\theta}\right\|\nonumber\\
    = & \sqrt{n}\left\|\proxi{f_{n}}{C_{n}}(\hat{\theta}_{init}) - \proxi{f_{n}}{C_{n}}(\hat{\theta})\right\| + \sqrt{n}\left\|\proxi{f_{n}}{C_{n}}(\hat{\theta}) - \hat{\theta}\right\|\nonumber\\
    = & \sqrt{n}\left\|\proxi{f_{n}}{C_{n}}(\hat{\theta}_{init}) - \proxi{f_{n}}{C_{n}}(\hat{\theta})\right\| + o_{P}(1). \label{prox_and_op}
  \end{align}
  We will focus our remaining efforts on the first term in \eqref{prox_and_op}. Assume that $\sqrt{n} \|\hat{\theta}_{init} - \hat{\theta}\| < M$ for some $M$, since its complement can be made to have negligible probability. By assumption $(ii)$ we can further assume that $\prox_{f_{n}}^{C_{n}}$ is single-valued. Finally, Lipschitz continuity of the scaled prox and assumption $(iii)$ allow us to assume that $f_{n}$ is strongly convex with modulus $m$ on some set containing both $\prox_{f_{n}}^{C_{n}}(\hat{\theta}_{init})$ and $\prox_{f_{n}}^{C_{n}}(\hat{\theta})$. From the definition of the proximal operator we have
  \begin{equation}\label{eq:proxcontain}
    \hat{\theta}_{init} \in \proxi{f_{n}}{C_{n}}(\hat{\theta}_{init}) + C_{n}^{-1} \partial f_{n}\left( \proxi{f_{n}}{C_{n}}(\hat{\theta}_{init})\right) \quad \text{and} \quad \hat{\theta} \in \proxi{f_{n}}{C_{n}}(\hat{\theta}) + C_{n}^{-1} \partial f_{n}\left(\proxi{f_{n}}{C_{n}}(\hat{\theta})\right).
  \end{equation}
We also have
  \begin{align}\label{eq:another_equation}
    & \lambda_{max}(C_{n}) \left\| \proxi{f_{n}}{C_{n}}(\hat{\theta}_{init}) - \proxi{f_{n}}{C_{n}}(\hat{\theta}) \right\| \left\| \hat{\theta}_{init} - \hat{\theta} \right\| \nonumber\\
     \geq & \left\| \proxi{f_{n}}{C_{n}}(\hat{\theta}_{init}) - \proxi{f_{n}}{C_{n}}(\hat{\theta}) \right\|_{C_{n}} \left\| \hat{\theta}_{init} - \hat{\theta} \right\|_{C_{n}},
  \end{align}
  so from \eqref{eq:proxcontain}, there are $u \in \partial f_{n}\left(\prox_{f_{n}}^{C_{n}}(\hat{\theta}_{init})\right)$ and $v \in \partial f_{n}\left(\prox_{f_{n}}^{C_{n}}(\hat{\theta})\right)$ such that \eqref{eq:another_equation} is bounded below as
  $$\geq \left\| \proxi{f_{n}}{C_{n}}(\hat{\theta}_{init}) - \proxi{f_{n}}{C_{n}}(\hat{\theta}) \right\|_{C_{n}} \left\| \proxi{f_{n}}{C_{n}}(\hat{\theta}_{init}) - C_{n}^{-1} u - \left(\proxi{f_{n}}{C_{n}}(\hat{\theta}_{init}) - C_{n}^{-1} v\right) \right\|_{C_{n}}$$
  Applying the Cauchy-Schwarz inequality
  \begin{align*}
    \geq & \left(\proxi{f_{n}}{C_{n}}(\hat{\theta}_{init}) - \proxi{f_{n}}{C_{n}}(\hat{\theta})\right)^{T} C_{n} \left( \proxi{f_{n}}{C_{n}}(\hat{\theta}_{init}) - \proxi{f_{n}}{C_{n}}(\hat{\theta}) + C_{n}^{-1} \left( u - v \right)\right)\\
    = & \left\|\proxi{f_{n}}{C_{n}}(\hat{\theta}_{init}) - \proxi{f_{n}}{C_{n}}(\hat{\theta})\right\|_{C_{n}}^{2} + \left(\proxi{f_{n}}{C_{n}}(\hat{\theta}_{init}) - \proxi{f_{n}}{C_{n}}(\hat{\theta})\right)^{T} \left(u - v\right)
  \end{align*} 
  Invoking strong convexity of $f_{n}$, we have that $\partial f_{n}$ is strongly monotone for some constant $m$, so we continue as follows.
  \begin{align*}
    \geq & \left\|\proxi{f_{n}}{C_{n}}(\hat{\theta}_{init}) - \proxi{f_{n}}{C_{n}}(\hat{\theta}) \right\|_{C_{n}}^{2} + m \left\| \proxi{f_{n}}{C_{n}}(\hat{\theta}_{init}) - \proxi{f_{n}}{C_{n}}(\hat{\theta}) \right\|^{2}\\
    \geq & \left(\lambda_{min}(C_{n}) + m\right) \left\|\proxi{f_{n}}{C_{n}}(\hat{\theta}_{init}) - \proxi{f_{n}}{C_{n}}(\hat{\theta}) \right\|^{2}
  \end{align*}
  From this chain of inequalities we conclude
    $$\frac{\lambda_{max}(C_{n})}{\lambda_{min}(C_{n}) + m} \left\|\hat{\theta}_{init} - \hat{\theta}\right\| \geq \left\|\proxi{f_{n}}{C_{n}}(\hat{\theta}_{init}) - \proxi{f_{n}}{C_{n}}(\hat{\theta})\right\|.$$
  Therefore if $\lambda_{max}(C_{n}) \to^{P} 0$, since $\sqrt{n}\left\|\hat{\theta}_{init} - \hat{\theta}\right\|$ is bounded in probability, we have $$\sqrt{n}\left\|\proxi{f_{n}}{C_{n}}(\hat{\theta}_{init}) - \proxi{f_{n}}{C_{n}}(\hat{\theta})\right\| = o_{P}(1).$$
  We conclude that
  $$\sqrt{n} \left\|\hat{\theta}_{ose} - \hat{\theta}\right\| = o_{P}(1),$$
  which completes the proof.
\endproof

Similar to Theorem \ref{thm:proxgrad}, we note that $\hat{\theta}$ need not be a global minimizer of $f_{n}$ to apply Theorem \ref{thm:proxdescent}, but instead satisfy the approximate fixed-point condition $\prox_{f_{n}}^{C_{n}}(\hat{\theta}) - \hat{\theta} = o_{P}(n^{-1/2})$. 

We also include a proposition which permits the application of Theorem \ref{thm:proxdescent} in the setting of parametric estimation from i.i.d. observations introduced in subsection \ref{subsec:one-step}. 
\begin{proposition}\label{prop:assump2}
  %The assumption of this proof follow that of Ferguson's Theorem 18
  Let $X_{1},...,X_{n} \sim^{i.i.d.} p_{\theta^{*}}$ and $f_{n} =  -\frac{1}{n} \sum_{i=1}^{n} \log p_{\theta}(X_{i})$ Assume the following.
  \begin{enumerate}[(i)]
    \item $\log p_{\theta} \in C^{2}$ for each $X$. Further, partial derivatives may be passed under the integral in $\int p_{\theta}(x) \, dx$.
    \item There exists a function $K(x)$ such that $E \left[K(X_{1})\right] < \infty$ and each component of $\nabla^{2} \log p_{\theta_{0}}$ is bounded in absolute value by $K$ uniformly in some neighborhood of $\theta_{0}$. %from theorem 18 in Ferguson
    \item The matrix $E[\nabla^{2} \log p_{\theta_{0}}]$ exists and is negative definite (so that the Fisher information $I_{\theta_{0}}$ exists and is nongenerate). Further, the Fisher information is continuous at $\theta_{0}$.
  \end{enumerate}
  Then Theorem \ref{thm:proxdescent} holds with only the first two assumptions of that theorem. That is, if
\begin{enumerate}[(i')]
\item $\lambda_{max}(C_{n}) = o_{P}(1)$ 
\item There is a constant $L$ such that for each constant $M$, the function $\prox_{f_{n}}^{C_{n}}$ is L-Lipschitz continuous on $\left\{\theta: \|\theta - \hat{\theta}\| \leq M/\sqrt{n}\right\}$ with high probability.
\end{enumerate}
  Then for any $\sqrt{n}$-consistent estimator $\hat{\theta}_{init}$ of $\hat{\theta}$, $\hat{\theta}$ from \eqref{eq:estproxd} and $\hat{\theta}_{ose}$ from \eqref{eq:oseproxd} are asymptotically equivalent.
\end{proposition}

\proof{}
For brevity, we denote $\log p_{\theta}$ by $\ell_{\theta}$. We first show that there exists $m > 0$ such that for all $\|h\| \leq B$, it holds with high probability that $f_{n} := - \sum_{i=1}^{n} \ell_{\theta_{0} + h/\sqrt{n}}(X_{i})$ is strongly convex as a function of $h$ with parameter $m$. Take $m$ to be such that $\lambda_{min}(I_{\theta_{0}}) > m > 0$.

  Assumption $(i)$ gives the well known result that $I_{\theta} = - E[\nabla^{2} \ell_{\theta}(X_{1})]$. By Taylor's theorem,
  $$\sum_{i=1}^{n} \ell_{\theta_{0} + h/\sqrt{n}}(X_{i}) = \sum_{i=1}^{n} \ell_{\theta_{0}}(X_{i}) + h^{T}\left(\frac{1}{\sqrt{n}} \sum_{i=1}^{n} \nabla \ell_{\theta_{0}}(X_{i})\right) + \frac{1}{2} h^{T} \left(\frac{1}{n} \sum_{i=1}^{n} \nabla^{2} \ell_{\theta_{0} + t_{i} h/\sqrt{n}}(X_{i}) \right)^{T}h $$
  for some $t_{i} \in (0,1)$. We note that the $t_{i}$ are random because they depend on $X_{i}$.

Our analysis will focus on the quadratic term. We apply the uniform law of large numbers \cite[Theorem 16a]{ferguson} and the continuity of $I_{\theta}$ to see that
\begin{align*}
  & \sup_{\|h\| \leq B} \left\|\frac{1}{n} \sum_{i=1}^{n} \nabla^{2} \ell_{\theta_{0} + t_{i} h/\sqrt{n}}(X_{i}) -E_{\theta^{*}} \left[ \nabla^{2} \ell_{\theta_{0}}(X_{1})\right] \right\|\\
  \leq & \sup_{\|h\| \leq B} \left\|\frac{1}{n} \sum_{i=1}^{n} \nabla^{2} \ell_{\theta_{0} + h/\sqrt{n}}(X_{i}) -E_{\theta^{*}} \left[ \nabla^{2} \ell_{\theta_{0}}(X_{1})\right] \right\|\\
\leq 
  & \sup_{\|h\| \leq B} \left\|\frac{1}{n} \sum_{i=1}^{n} \nabla^{2} \ell_{\theta_{0} + h/\sqrt{n}}(X_{i}) - E_{\theta^{*}} \left[ \nabla^{2} \ell_{\theta_{0} + h/\sqrt{n}}(X_{1})\right] \right\| \\
  & + \sup_{\|h\| \leq B} \left\|
  E_{\theta^{*}} \left[ \nabla^{2} \ell_{\theta_{0} + h/\sqrt{n}}(X_{1})\right] - 
  E_{\theta^{*}} \left[ \nabla^{2} \ell_{\theta_{0}}(X_{1})\right]
 \right\|\\
  = & \sup_{\|h\| \leq B} \left\|\frac{1}{n} \sum_{i=1}^{n} \nabla^{2} \ell_{\theta_{0} + h/\sqrt{n}}(X_{i}) + - E_{\theta^{*}} \left[ \nabla^{2} \ell_{\theta_{0} + h/\sqrt{n}}(X_{1})\right] \right\| \\
  & + \sup_{\|h\| \leq B} \left\|
  I_{\theta_{0} + h/\sqrt{n}} - 
  I_{\theta_{0}} \right\|.\\
\end{align*}
  The first inequality in the above display follows from $t_{i} \in (0,1)$. Since both of the terms in the last equation converge to zero in probability, the Hessian term in the Taylor expansion of $\ell_{\theta_{0} + h/\sqrt{n}}$ converges in probability to $-I_{\theta_{0}}$. Thus for each $\epsilon > 0$, we can choose $N$ such that $n \geq N$ implies
  $$\frac{1}{n} \sum_{i=1}^{n} \nabla^{2} \ell_{\theta_{0} + t_{i} h/\sqrt{n}}(X_{i}) \preceq - I_{\theta_{0}} + \epsilon I$$
  with arbitrarily high probability. Thus $f_{n}$ is strongly convex with high probability for $m > 0$.

  We have shown that, for each constant $M$, the function $f_{n}$ is strongly convex with high probability on $\{\theta: \norm{\theta - \hat{\theta}} \leq M/\sqrt{n}\}$. This is assumption $(iii)$ in Theorem \ref{thm:proxdescent}, so the result is proven.
\endproof

\section{The Moreau Envelope as a Statistical Smoother.}\label{Moreau}

In this section we interpret our contributions in the context of smoothing irregularities in a statistical objective. Theorems \ref{thm:proxgrad} and \ref{thm:proxdescent} both provide equivalences between OSE and MLE estimators, but additionally allow an important connection to infimal convolution and the scaled Moreau envelope as a statistical smoother.

For two extended real-valued functions $f_{1}: \R^{d} \to \R \cup \{\infty\}$ and $f_{2}: \R^{d} \to \R \cup \{\infty\} $, the infimal convolution $f_{1} \# f_{2}: \R^{d} \to \R \cup \{\infty\}$ is
$$(f_{1} \# f_{2})(x) = \inf_{x_{1} + x_{2} = x} \left\{ f_{1}(x_{1}) + f_{2}(x_{2}) \right\}$$
As long as the infimum in the definition of $(f_{1} \# f_{2})$ is attained, the epigraph of the inf-convolution is the Minkowski sum of the epigraphs of $f_{1}$ and $f_{2}$ \cite{VarAn}. This interpretation allows the infimal convolution to be used as a smoothing operation, where $f_{1} \# f_{2}$ gives a smoothed version of $f_{1}$ for $f_{2}$ chosen with an epigraph satisfying certain regularity properties. See any of \cite{Burke2013Epi, Xu2014Smoothing, Burke2016Epi} for recent development of inf-convolution as a smoother, \cite{VarAn} for a more general overview of the properties of inf-convolution, and \cite{BeTe12smoothing} for more on Moreau envelope as a smoother, as well as its connection to the proximal operator.

The Moreau envelope is a special case of inf-convolution, where a function is inf-convolved against a squared $\ell_{2}$ norm. For $f: \R^{d} \to \R \cup \{\infty\}$ and positive definite matrix $C \in \R^{d \times d}$, the scaled Moreau envelope $e_{C} f: \R^{d} \to \R$ is
$$e_{C} f(x) = \inf_{w \in \R^{d}}\left\{ f(w) + \frac{1}{2} \|x - w\|_{C}^{2} \right\}.$$
The (unscaled) Moreau envelope is the scaled Moreau envelope with $C = (1/\lambda) I$ for some positive scalar $\lambda$. The scaled Moreau envelope is closely linked to the scaled prox, because for any $w^{*} \in \prox_{f}^{C}(x)$,
$$e_{C} f(x) = f(w^{*}) + \frac{1}{2} \| x - w^{*}\|_{C}^{2}.$$

Despite the tendency of minimization to destroy smoothness, the Moreau envelope is smooth for convex functions $f$ \cite{VarAn}. In the nonconvex case, smoothness of the Moreau envelope at a point $x^{*}$ requires inner continuity\footnote{In the set-convergence sense of variational analysis. See \cite{VarAn}. If $\prox_{f}^{C}$ is single-valued in a neighborhood of $x^{*}$ then inner-continuity reduces to continuity.} of the scaled prox at a point in $\{x^{*}\} \times \prox_{f}^{C}(x^{*})$. We formalize this result in the following proposition.
\begin{proposition}{\rm (differentiability of $e_{C}f$)}\label{prop:diffeCf}
Let $f:\R^d\to\R\cup\{+\infty\}$ be an extended real-valued function.  For a positive definite matrix $C \succ0$, the Moreau envelope $e_{C}f$ is differentiable at $x^*$ if there exists $w^*\in\prox_{f}^{C}(x^*)$ such that for all sequences $u^\nu\to 0$ there exists
  $w^\nu\in\prox_{{f}}^{{C}}(x^*+u^\nu)$ with $w^\nu \to w^*$. In this case  $\nabla e_{C}f(x^*) = C(x^* - w^*)$. 
\end{proposition}

\proof{} See appendix. \endproof
 
 \smallskip

\begin{figure}[ht]
	\begin{center}
		\includegraphics[scale=.6]{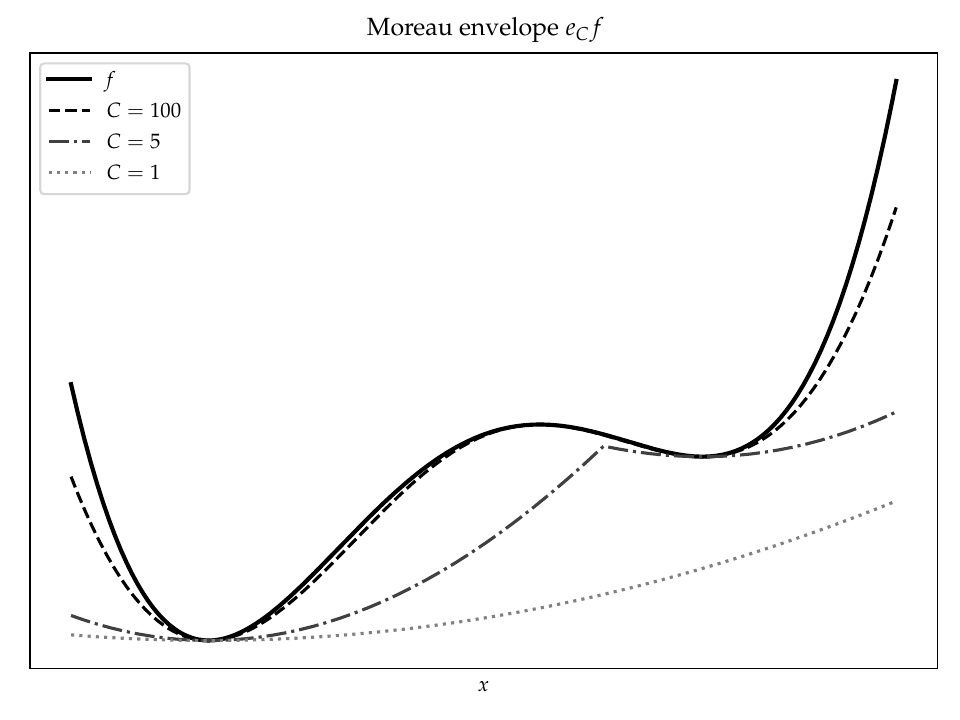}
		\caption{A collection of (unscaled) Moreau envelopes for a non-convex function}
		\label{fig:Moreau}
	\end{center}
\end{figure}

It is clear from Proposition~\ref{prop:diffeCf} that the objective function $e_Cf$ is a smoothed version of the function $f$ when $f$ satisfies the required assumptions. Moreover, $e_{C}f$ is finite over $\R^{d}$ when $f$ is proper and bounded from below, even though $f$ itself might take the value $\infty$. We also note that $e_{C} f$ preserves global minimizers of $f$; if $x^*$ is a minimizer of $f$ then the following chains of inequalities hold for all $x$ and $w$.
$$f(x^*)\leq f(w)+\frac{1}{2}\|x-w\|_{C}^2$$
$$e_C f(x^*) \leq f(x^{*}) + \frac{1}{2} \| x^{*} - x^{*}\|^{2} = f(x^{*}).$$
Taking the infimum over $w$ and $x$ implies that the sets of minimizers of $f$ and $e_{C}f$ coincide. The smoothing properties of the Moreau envelope are depicted in Figure \ref{fig:Moreau}, which also illustrates the coincidence of minimizers.
% We can verify the smoothness of the Moreau Envelope by means of convex analysis.

When the Moreau envelope is differentiable at a point $x$, Proposition \ref{prop:diffeCf} gives
\[ \nabla e_{C} f(x)  =  C\left(x - \proxi{f}{C}(x)\right)\,\Rightarrow\,
 \proxi{f}{C}(x)  =  x - C^{-1}  \nabla e_{C} f(x) \]
We conclude that the proximal operator is scaled gradient descent applied to the Moreau envelope. This perspective allows us to interpret Theorem \ref{thm:proxdescent} in the context of smoothing the likelihood function with the Moreau envelope.

\begin{corollary} \label{cor:moreausmooth}
  Under the assumptions of Theorem \ref{thm:proxdescent}, the one-step estimator 
  $$\hat{\theta}_{ose}= \hat{\theta}_{init} - C^{-1}_{n} \nabla e_{C_{n}} f_{n}(\hat{\theta}_{init}),$$
  which is a single iteration of scaled gradient descent applied to the Moreau envelope $e_{C_{n}}$, is asymptotically equivalent to the estimator $\hat{\theta}$ which minimizes $f_{n}$
\end{corollary}

Interpreted in the context of the Moreau envelope, Corollary \ref{cor:moreausmooth} provides a recipe for smoothing a negative log-likelihood, while retaining the attractive theoretical properties of the maximum likelihood estimator. Thus, the Moreau envelope can be used to smooth a statistical objective, with the goal of removing irrelevant local artefacts while preserving important global structure. Other efforts in this area include the continuation method, also called graduated optimization, in which an objective is convolved against a smooth (usually Gaussian) kernel in order to impart additional smoothness in the objective function \cite{Hazan2016}, \cite{Mobahi}. The continuation method also appears independently in the statistics literature as Maximum Smoothed Likelihood Estimation \cite{Ionides, Eggermont2001}. The main limitation to the practical application of the continuation method is the difficulty associated with computing gradient and curvature information for the smoothed function. Smoothing objectives with the Moreau envelope, on the other hand, has the attractive property that descent steps are evaluations of the scaled proximal operator. Many techniques exist to compute the scaled proximal operator; for examples see \cite{FriedlanderGoh, glmnet, leesunsaunders}.

\section{Experiments \& Examples.}\label{Numerics}

\subsection{Estimation with the Cauchy Distribution.}
In this section we discuss our results in the context of estimating the location parameter of a Cauchy distributed random variable. We begin by formulating a penalized objective using a MAP (Maximum-a-Posteriori) estimator, to which we apply scaled proximal gradient descent. We then alter the problem to remove the prior information, and apply scaled proximal descent and Theorem \ref{thm:proxdescent} to the MLE problem.

\begin{figure}[ht]
\begin{minipage}{.5\linewidth}
\centering
\subfloat[]{\label{proxgrad1}\includegraphics[scale=.5]{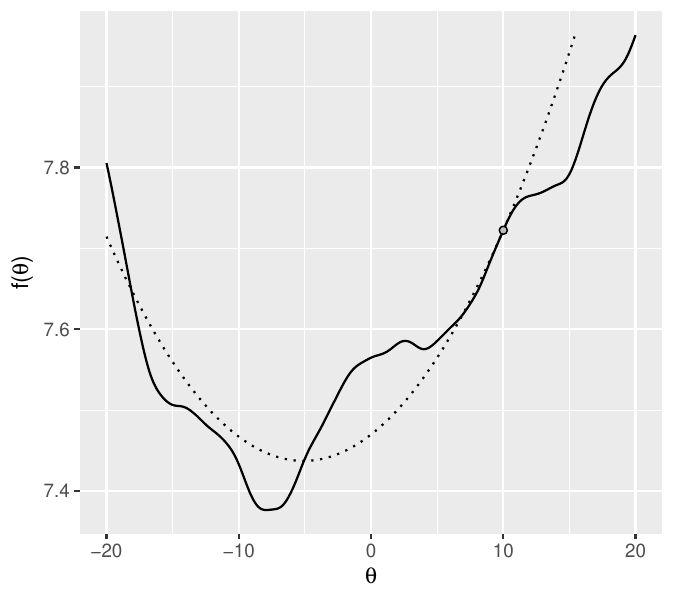}}
\end{minipage}
  \hspace{-.2cm}
\begin{minipage}{.5\linewidth}
\centering
  \subfloat[]{\label{proxgrad2}\includegraphics[scale=.5]{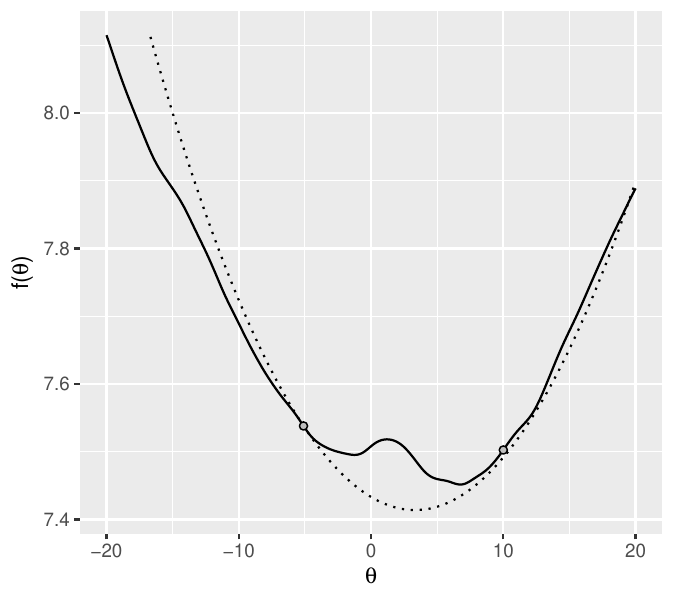}}
\end{minipage}\\
\begin{minipage}{.5\linewidth}
\centering
\subfloat[]{\label{proxgrad3}\includegraphics[scale=.5]{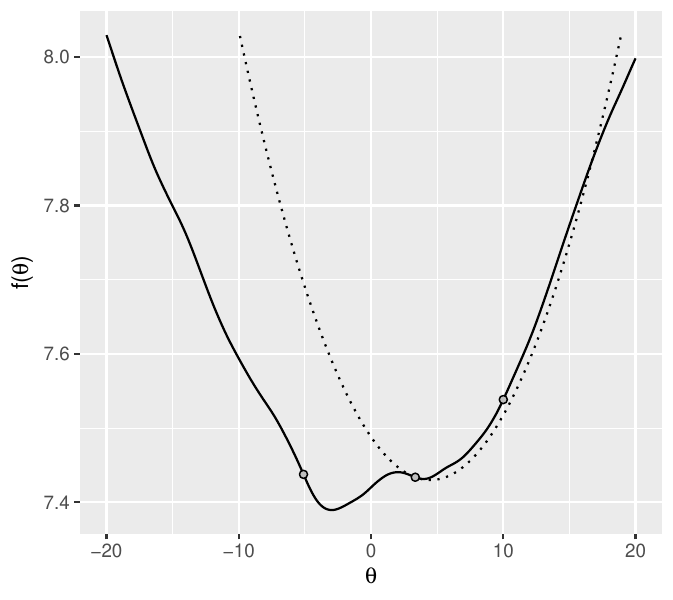}}
\end{minipage}%
\begin{minipage}{.5\linewidth}
\centering
\subfloat[]{\label{proxgrad4}\includegraphics[scale=.5]{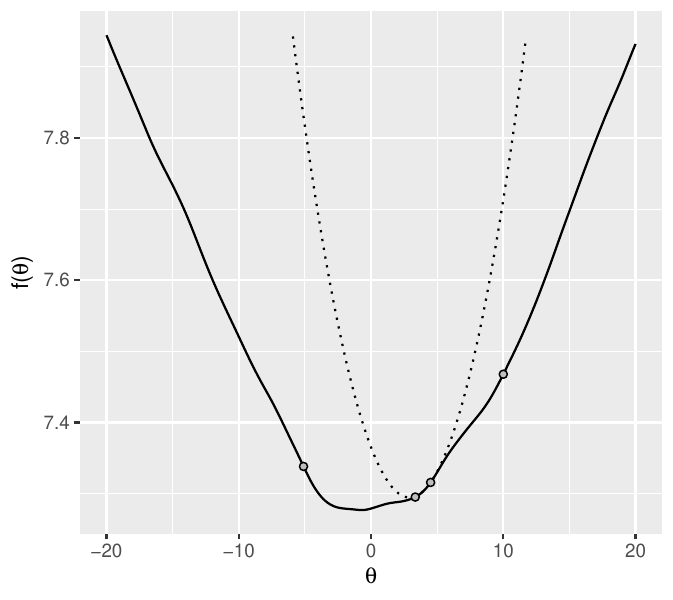}}
\end{minipage}
  \caption{A sequence of iterates of proximal gradient descent for increasing $n$ and $C_{n} \to 1/\sigma_{0} = 80$. The solid line gives the objective function, which changes as $n$ increases. The dashed line gives the function that is minimized to determine the next iterate. Choosing $C_{n}$ small in early iterates has a smoothing effect on the sequence.}
\label{fig:proxgrad}
\end{figure}

Let $X_{1},..., X_{n} \sim^{i.i.d.} \mathop{Cauchy}(\theta_{0}, \sigma_{0})$, where $\theta_{0}$ is the location parameter and $\sigma_{0}$ the scale parameter. Assume that $\sigma_{0}$ is known and we wish to estimate $\theta_{0}$. Moreover, assume that we incorporate a regularizer through a Laplacian prior on $\theta_{0}$, where $\sigma_{0} \sim \mathop{Laplace}(0, \gamma)$. The sample and prior information can be combined using a Maximum-a-Posteriori formulation, where the objective takes the form
$$\min_{\theta} \; - \frac{1}{n} \sum_{i=1}^{n} \left\{\log\left(1 + \left(\frac{X_{i} - \theta}{\sigma_{0}}\right)^{2}\right)\right\} + \frac{1}{n \gamma}\abs{\theta}.$$
This problem is nonconvex and has multiple local minima, making its computation through iterative methods difficult. Moreover, the sample mean is not consistent for the location parameter. Despite these challenges, the Cauchy distribution yields a negative log-likelihood $g_{n}$ which satisfies the conditions of Proposition \ref{prop:assump2}. Hence $\theta_{0}$ can be estimated using a one-step scaled proximal gradient estimator, where the scaling $C_{n} \to \frac{1}{2\sigma_{0}}$, the Fisher information of the Cauchy distribution at $\theta = 0$. We will use this example to illustrate the implications of Theorem \ref{thm:proxgrad} for nonconvex problems generally. In the remainder we fix $\sigma_{0} = 20$ and $\gamma = 1000$. 

Recall that iterations of scaled proximal gradient descent can also be written
$$\theta_{k+1} = \argmin_{\theta} \; \left\{ g_{n}(\theta_{k}) + \nabla g_{n}(\theta_{k})^{T}(\theta - \theta_{k}) + h_{n}(\theta) + \left\|\theta - \theta_{k}\right\|_{C_{n}}^{2} \right \}.$$
Hence iterates can be viewed as solving a penalized version of the original problem, where $g_{n}$ is linearized around the current iterate. Because this example is univariate, scaled proximal gradient descent is equivalent to the unscaled version. Figure \ref{fig:proxgrad} gives a sequence of iterates of scaled proximal gradient descent, as the scaling increases from $1/400$ to $1/2\sigma_{0}$ and the sample size $n$ increases from $100$ to $1000$.

\begin{figure}[ht]
	\begin{center}
		\includegraphics[scale=.9]{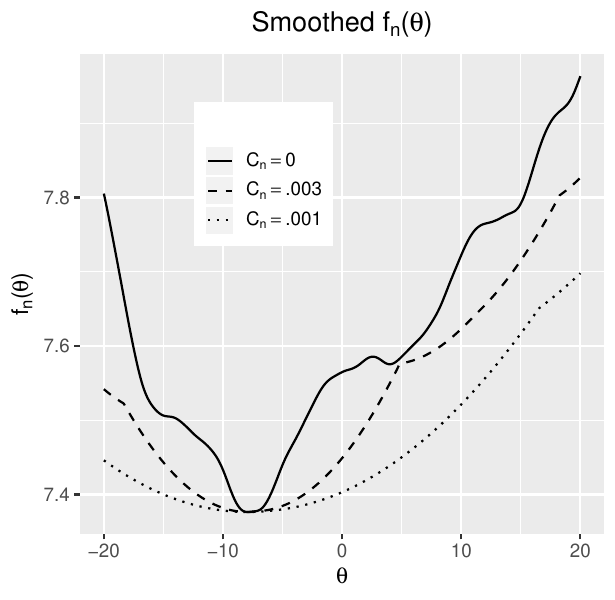}
		\caption{Negative log-likelhood and Moreau envelopes for the location parameter in a Cauchy distributed random variable. In this example, $\theta_{0} = 0$, $\sigma_{0} = 20$, and $n=100$.}
		\label{fig:mleCauchy}
	\end{center}
\end{figure}

%\lambda is C_{n}^{-1} in the connection between scaled and unscaled
If we omit the Laplacian prior, then we have a maximum likelihood objective without regularization. In this case, the smoothing discussion in the previous section applies, so that scaled proximal descent applied to this objective is scaled gradient descent applied to the Moreau envelope of the negative log-likelihood. Interpreted in this context, Theorem \ref{thm:proxdescent} provides theoretical justification for smoothing the Cauchy distribution's negative log-likelihood with the Moreau envelope. Figure \ref{fig:mleCauchy} illustrates this smoothing for a certain sample.

%For this particular sample, the function $\ell$ exhibits a local minimum on the left tail of the distribution, and the global minimum near the origin.  In particular, if we are interested in solving the critical points of $\ell$ by solving Equation~(\ref{eq:rootcau}), there is no guarantee that we will get the optimal points for this problem.  As suggested in \cite[Ex.5.5]{vanderVaart}, in order to obtain the maximum likelihood estimator, we should proceed using an algorithmic approach. Despite the fact that the sample mean is not consistent, the sequence of maximum likelihood estimators is consistent. Moreover, from the Cram\'er-Rao asymptotic consistency conditions, it is easy to see that the Fisher information for a sample of size $n$ is given by $\frac{2}{n}$, and then, the MLE $\hat{\theta}$ for this problem is $\sqrt{n}$-asymptotically normal \cite[Ch.16]{Das08asymp}.
%
%Finally, to see the full capabilities of the application of the one step estimator, we collect all the information that we have available: the MLE for a localized Cauchy distribution function is consistent and $\sqrt{n}$-asymptotically normal, and is the solution of the minimization problem given by equation~(\ref{mle:Cauchy}), with a nonconvex objective function.  Our proposed strategy of a one-step estimator, using the prox operator described on Section~\S\ref{MLE} provides a numerical procedure to solve this problem to an asymptotically equivalent solution, in virtue of Theorem~\ref{thm:mle}
%
\subsection{First Order Methods Are Not Enough.}
The asymptotic equivalence of the scaled proximal gradient one-step estimator and the MLE prompts the following question: can a one-step estimator without the scaling, such as gradient descent, also be shown to be asymptotically equivalent to the MLE? In this section we will provide a counterexample demonstrating that without curvature information the gradient one-step estimator cannot overcome bias introduced in the distribution of $\hat{\theta}_{init}$.

Assume we want to estimate the mean of a bivariate normal distribution with known variance. Fix $\mu \in \R^2$ and let $X_1,.., X_{n} \sim^{i.i.d.} N(\mu, \Sigma)$, where $$\Sigma = \diag\left(\left( \begin{array}{cc} \sigma_{1}^{2}, \sigma_{2}^{2} \end{array} \right)^{\top} \right).$$
Assume without loss of generality that $\sigma_{1}^{2} > \sigma_{2}^{2}$. We denote by $\hat{\theta}_{init}$ the initial estimate and $\hat{\theta}$ the maximum likelihood estimator of the model. The negative log-likelihood function is quadratic with Hessian $\Sigma^{-1}$. Let $\alpha$ be the step length of gradient descent--to be specified momentarily. Choose $\hat{\theta}_{init} \sim U\left((\mu_{1} - \frac{1}{\sqrt{n}}, \mu_{1}) \times (\mu_{2} - \frac{1}{\sqrt{n}}, \mu_{2})\right)$ independently of $X_{1},...,X_{n}$. We have, up to a constant factor
$$f_{n}(\theta) = \frac{1}{n}\sum_{i=1}^{n} \frac{1}{2}(\theta - X_{i})^{\top} \Sigma^{-1} (\theta - X_{i}).$$
Hence
$$\hat{\theta} = \hat{\theta}_{init} - \alpha \nabla f_{n}(\hat{\theta}_{init}) = \hat{\theta}_{init} -\alpha \sigma^{-2}_{1} (\hat{\theta}_{init} - \bar{X}).$$

Consider $P(\hat{\theta}_{1} > \mu_{1})$, the probability that $\hat{\theta}$ exceeds $\mu$ in its first coordinate (where we denote the first coordinate with a subscript).
$$P(\hat{\theta}_{1} > \mu_{1}) = P(\hat{\theta}_{init, 1} - \alpha \Sigma^{-1} (\hat{\theta}_{init, 1} - \bar{X}_{1}) > \mu_{1})$$
$$=P((1-\alpha\sigma_{1}^{-2})(\hat{\theta}_{init, 1}-\mu_{1}) + \alpha \sigma_{1}^{-2}(\bar{X}_{1} - \mu_{1}) > 0)$$
$$=P\left(\bar{X}_{1} - \mu_{1} > \frac{(1-\alpha \sigma_{1}^{-2})(\mu_{1} - \hat{\theta}_{init, 1})}{\alpha \sigma_{1}^{-2}}\right)$$
We have $\sqrt{n} (\bar{X}_{1} - \mu_{1}) = \sigma_{1} Z$, where $Z$ is a standard normal. Further, $\sqrt{n} (\mu_{1} - \hat{\theta}_{init, 1}) = U$, where $U$ is uniform on $[0,1]$. We multiply both sides of the inequality above by $\sqrt{n}$ and rewrite as follows.
$$=P\left(Z > \frac{(1-\alpha \sigma_{1}^{-2}) U}{\alpha \sigma_{1}^{-1}}\right)$$
$$=P\left(Z > \left(\frac{\sigma_{1}}{\alpha} - \sigma_{1}^{-1}\right) U\right).$$

Note that if we can make $\sigma_{1}/\alpha - \sigma_{1}^{-1} = M$ for some fixed large $M$, then integration of the above expression shows that
$$P(\hat{\theta}_1 > \mu_{1}) = 1-\Phi(M) + (1-\exp(-M^{2}/2)/\sqrt{2\pi})/M.$$
We remark that this probability is independent of $n$, and can be made arbitrarily small for $M$ chosen large.

We will consider two common choices of step length $\alpha$, fixed step length and exact line search. The objective function $f$ has Lipschitz constant $\sigma_{2}^{-2}$. According to \cite{linearnonlinear}, the optimal fixed step length in gradient descent is the reciprocal the objective's Lipschitz constant; $\sigma_{2}^{2}$ in our example. Therefore
$$\frac{\sigma_{1}}{\alpha} - \sigma_{1}^{-1} = \frac{\sigma_{1}}{\sigma_{2}^{2}} - \sigma_{1}^{-1}.$$
By fixing $\sigma_{1}$ and choosing $\sigma_{2}$ small, this expression can be made large, so that the probability $P(\hat{\theta}_{1} > \mu_{1})$ is arbitrarily large independent of $n$. Since the maximum likelihood estimator $\bar{X}$ has $P(\bar{X}_{1} > \mu_{1})$, $\hat{\theta}_{1}$ is not asymptotically equivalent to the MLE.

\begin{figure}[ht]
\begin{minipage}{.5\linewidth}
\centering
\subfloat[]{\label{main:a}\includegraphics[scale=.4]{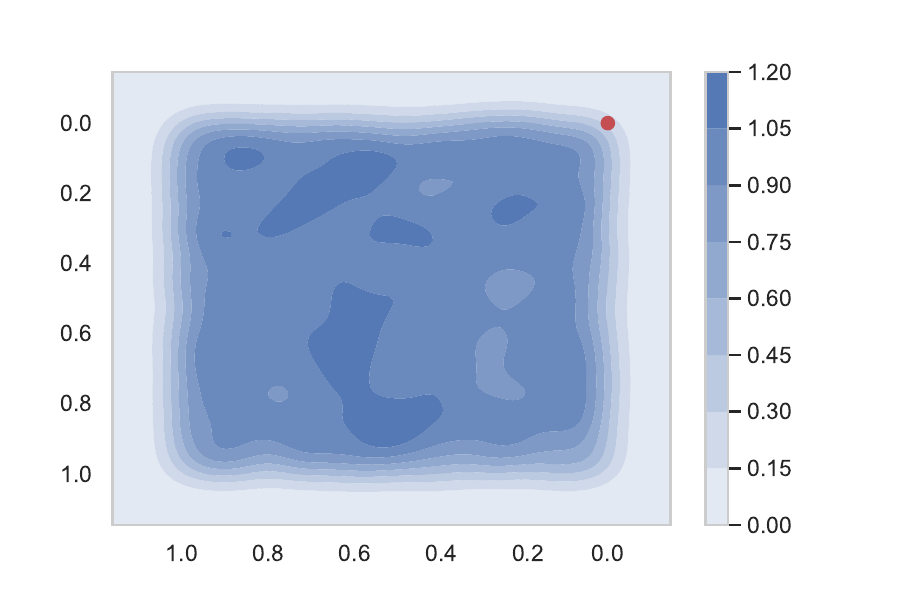}}
\end{minipage}
  \hspace{-.2cm}
\begin{minipage}{.5\linewidth}
\centering
  \subfloat[]{\label{main:b}\includegraphics[scale=.4]{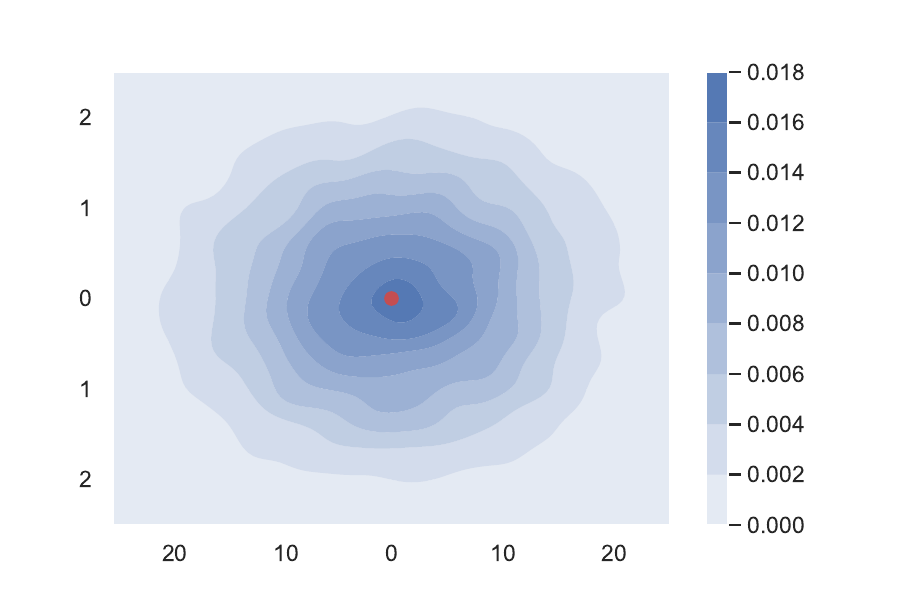}}
\end{minipage}\\
\begin{minipage}{.5\linewidth}
\centering
\subfloat[]{\label{main:c}\includegraphics[scale=.4]{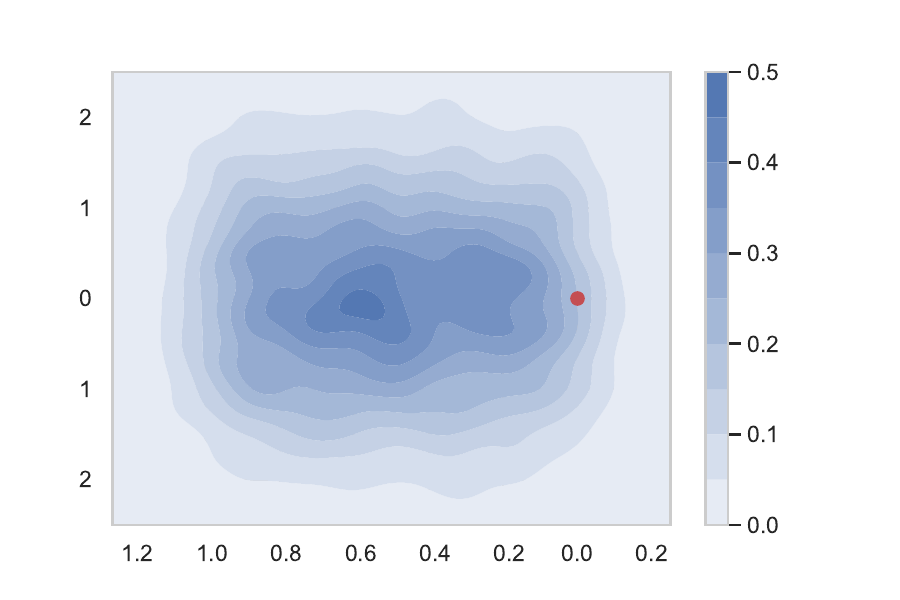}}
\end{minipage}%
\begin{minipage}{.5\linewidth}
\centering
\subfloat[]{\label{main:d}\includegraphics[scale=.4]{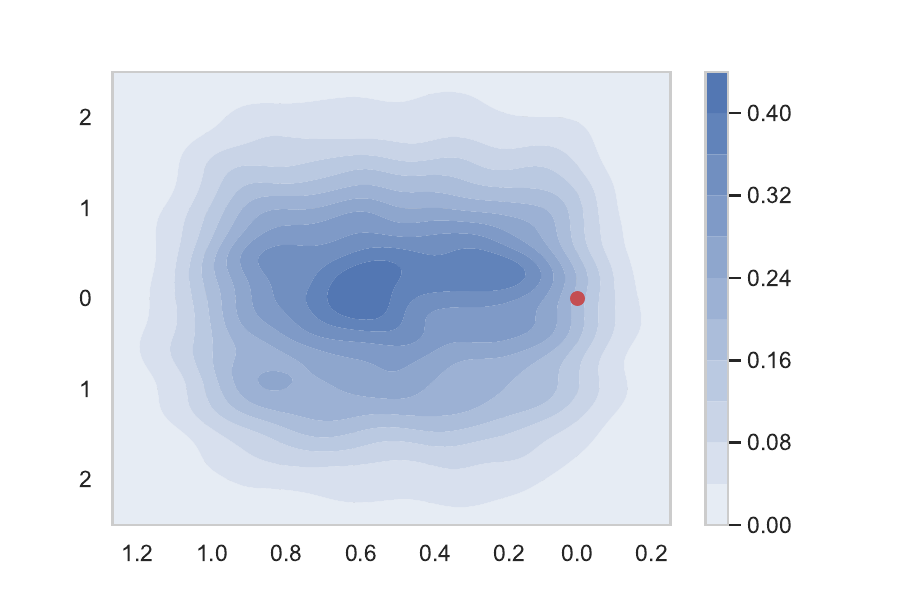}}
\end{minipage}
\caption{Kernel density estimates of the $\sqrt{n}$-normalized asymptotic distributions for (a) the starting point, (b) maximum likelihood estimator, (c) the one-step gradient descent estimator with optimal \emph{fixed} step length, and (d) the one-step step gradient descent estimator with optimal step length In each plot, the population mean of $X$ is given in red. Ten thousand samples used to construct each density estimate. The estimators in (b)-(d) were each constructed from samples of $X$ of size ten thousand. In this example, $\sigma_{1} = 10$ and $\sigma_{2} = 1$.}
\label{fig:fig1}
\end{figure}

Moreover, we can extend this result beyond the fixed step length to the setting of \emph{exact} step lengths. With a quadratic objective, the exact step length for steepest descent is (see \cite{linearnonlinear}),
$$\alpha := \frac{\nabla f_{n}(\theta)^{\top} \nabla f_{n}(\theta)}{\nabla f_{n}(\theta)^{\top} \Sigma^{-1} \nabla f_{n}(\theta)}.$$
In the example above this is
$$\alpha = \frac{(\bar{X} - \theta)^{\top} \Sigma^{-2} (\bar{X} - \theta)}{(\bar{X} - \theta)^{\top} \Sigma^{-3} (\bar{X} - \theta)}$$
%$$=\frac{\sigma_{1}^{-2}(\bar{X} - \theta_{1})^{2} + \sigma_{2}^{-2}(\bar{X} - \theta_{2})^{2}}{\sigma_{1}^{-3}(\bar{X}_{1} - \theta_{1})^{2} + \sigma_{2}^{-3}(\bar{X}_{2} - \theta_{2})^{2}}$$
$$=\frac{\sigma_{1}^{-2}(\bar{X} - \theta_{1})^{2}}{\sigma_{1}^{-3}(\bar{X}_{1} - \theta_{1})^{2} + \sigma_{2}^{-3}(\bar{X}_{2} - \theta_{2})^{2}} +  \frac{\sigma_{2}^{-2}(\bar{X} - \theta_{2})^{2}} {\sigma_{1}^{-3}(\bar{X}_{1} - \theta_{1})^{2} + \sigma_{2}^{-3}(\bar{X}_{2} - \theta_{2})^{2}}$$
$$\leq \sigma_{1} + \sigma_{2}.$$
Therefore
$$\frac{\sigma_{1}}{\alpha} - \sigma_{1}^{-1} \geq \frac{\sigma_{1}}{\sigma_{1} + \sigma_{2}} - \sigma_{1}^{-1}$$
so when $\sigma_{1}/(\sigma_{1} + \sigma_{2}) - \sigma_{1}^{-1} > 0$,
$$P(\hat{\theta}_{1} > \mu_{1}) =  P\left(Z > \left(\frac{\sigma_{1}}{\alpha} - \sigma_{1}^{-1}\right) U\right)$$
$$\leq P\left(Z > \left(\frac{\sigma_{1}}{\sigma_{1} + \sigma_{2}} - \sigma_{1}^{-1}\right) U\right).$$
For $M$ arbitrarily large, we can choose $\sigma_{1} = 1$ and $\sigma_{2}$ arbitrarily small so that $\frac{\sigma_{1}}{\sigma_{1} + \sigma_{2}} - \sigma_{1}^{-1} = M$. Therefore we can again make $P(\hat{\theta}_{1} > \mu_{1})$ arbitrarily small independently of $n$, which shows that again the one-step estimator $\hat{\theta}_{n}$ is not asymptotically equivalent to the MLE. Figure \ref{fig:fig1} gives the empirical distribution of the one-step gradient descent estimator for a number of samples $X$ and starting points $\hat{\theta}_{init}$, and for both the optimal fixed and optimal step lengths.

\subsection{Low-Rank Logistic Regression.} \label{lowrank}

We next consider an application in low-rank logistic regression, where we will apply our one-step estimation results for the composite model from section \ref{Composite}. By using one-step estimation  along with the stopping condition presented in Proposition \ref{stop_cond}, we construct an estimator which has the same large sample performance as the nuclear-norm penalized MLE, but which is achieved in a finite number of iterations of scaled proximal gradient descent.

%Is this \sqrt{n} consistent? Additional paragraph to clarify that the assumptions are met.

Consider the email-Eu-core data set \cite{eu-core}, formed by a collection of the time indices for emails sent within an academic department. We discretize the time into 49 segments, and we assume that the probability of an email being sent from individual $i$ to individual $j$ is constant throughout time and equal to $P_{i,j}$. Our goal is to estimate $P$, the $N \times N$ matrix of communication probabilities, where $N$ is the number of individuals in the department. We further assume that the logit of the communication probabilities, $\log\left(\frac{P}{1-P} \right)$, where $\log$ operates elementwise, is low-rank, reflecting that individuals may have similar communication pattern across all members of the department. We propose to estimate $P$ through penalized logistic regression. 

\begin{figure}[ht]
    \centering
    \begin{minipage}{0.5\textwidth}
        \centering
        \includegraphics[width=.95\textwidth]{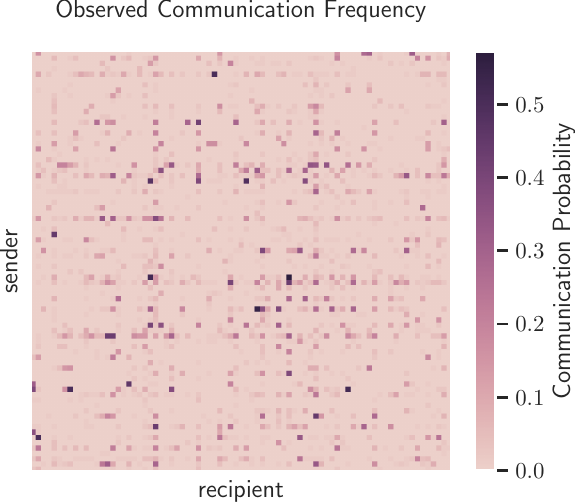} % second figure itself
	\caption*{(a)}
    \end{minipage}\hfill
    \begin{minipage}{0.5\textwidth}
        \centering
        \includegraphics[width=.95\textwidth]{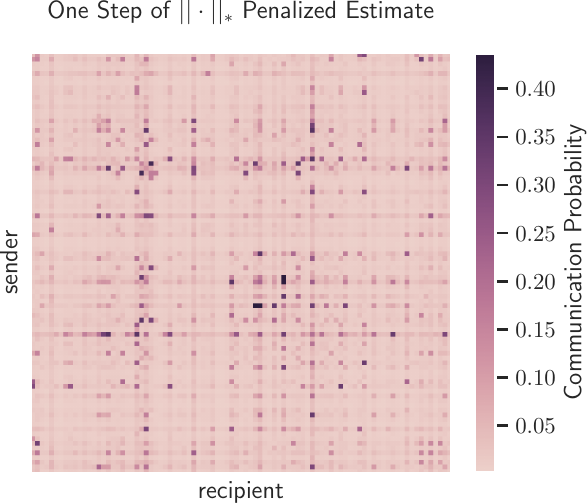} % first figure itself
	\caption*{(b)}
    \end{minipage}
    \caption{\label{fig:lowrank} (a) The empirical distribution for communication probability between each sender/receiver pair.
  (b) One step of proximal Newton descent applied to the nuclear-norm penalized logistic regression problem in section \ref{lowrank}.}
\end{figure}

\begin{equation}\label{logreg_model}
  \min_{\theta \in \mathbb{R}^{N \times N}} \sum_{i,j \in [N] \times [N]} \left\{\log(\exp(\theta_{i,j}) + 1) - \bar{X}_{i,j} \theta_{i,j} \right \} + \lambda \|\theta\|_{*}
\end{equation}
Above, $\|\cdot\|_{*}$ denotes the nuclear norm, and its inclusion as a penalty encourages low rank structure in $\theta$. The matrix $P$ can then be estimated by $\frac{\exp(\theta)}{1 + \exp(\theta)}$, where $\exp$ gives the elementwise exponential. We note that in general $\rank(\theta) \neq \rank(P)$, where $P = \exp(\theta)/(1+\exp(\theta))$, so the nuclear norm penalty on $\theta$ encourages underlying low-rank structure, though this low-rank structure is in $\theta$, the logit of $P$, instead of $P$.

We solve the optimization problem by implementing a proximal Newton algorithm, with a stopping criteria as in Proposition~\ref{stop_cond} to generate $\hat{\theta}_{init}$. We apply Proposition~\ref{stop_cond} and take $\hat{\theta}_{init}$ to be a point at which the step length of scaled proximal gradient iterations is less than $n^{-1/2}$. We take $C_{n}$ to be the Hessian at $\hat{\theta}_{init}$, so that $C_{n} \to I_{\theta^{*}}^{-1}$. The logistic regression model that we have formulated satisfies the assumptions of Proposition~\ref{prop:assump1} for $n$ large enough and penalty parameter $\lambda$ chosen appropriately because the regularity conditions (i)-(iii) are satisfied and the minimizer $\hat{\theta}$ of \eqref{logreg_model} converges to the true $\theta^{*}$ at a $\sqrt{n}$-rate \cite{penalized_mle}. This permits application of Theorem~\ref{thm:proxgrad}. Theorem \ref{thm:proxgrad} gives that the one-step estimator is asymptotically equivalent to the true minimizer of the nuclear-norm penalized logistic regression problem. Figure~\ref{fig:lowrank} displays our results, where the plot in (b) depicts the low rank structure of the solution obtained from applying one-step of proximal Newton descent.

\bibliographystyle{siamplain}
\bibliography{Draft1}

\begin{thebibliography}{10}

\bibitem{penalized_mle}
{\sc A.~Antoniadis, I.~Gijbels, and M.~Nikolova}, {\em Penalized likelihood
  regression for generalized linear models with non-quadratic penalties},
  Annals of the Institute of Statistical Mathematics, 63 (2011), pp.~585--615.

\bibitem{BassettDeride}
{\sc R.~Bassett and J.~Deride}, {\em Maximum a posteriori estimators as a limit
  of bayes estimators}, Mathematical Programming, 174 (2019), pp.~129--144.

\bibitem{BauCom11Cvx}
{\sc H.~H. Bauschke, P.~L. Combettes, et~al.}, {\em Convex analysis and
  monotone operator theory in Hilbert spaces}, vol.~408, Springer, 2011.

\bibitem{BeckBook}
{\sc A.~Beck}, {\em First-order methods in optimization}, vol.~25, SIAM, 2017.

\bibitem{BeTe09fista}
{\sc A.~Beck and M.~Teboulle}, {\em A fast iterative shrinkage-thresholding
  algorithm for linear inverse problems}, SIAM journal on imaging sciences, 2
  (2009), pp.~183--202.

\bibitem{BeTe12smoothing}
{\sc A.~Beck and M.~Teboulle}, {\em Smoothing and first order methods: A
  unified framework}, SIAM Journal on Optimization, 22 (2012), pp.~557--580.

\bibitem{sr1}
{\sc S.~Becker and J.~Fadili}, {\em A quasi-newton proximal splitting method},
  in Advances in Neural Information Processing Systems, 2012, pp.~2618--2626.

\bibitem{bickel1975one}
{\sc P.~J. Bickel}, {\em One-step huber estimates in the linear model}, Journal
  of the American Statistical Association, 70 (1975), pp.~428--434.

\bibitem{FiniteWilks}
{\sc S.~Boucheron and P.~Massart}, {\em A high-dimensional wilks phenomenon},
  Probability theory and related fields, 150 (2011), pp.~405--433.

\bibitem{Burke2013Epi}
{\sc J.~V. Burke and T.~Hoheisel}, {\em Epi-convergent smoothing with
  applications to convex composite functions}, SIAM Journal on Optimization, 23
  (2013), pp.~1457--1479.

\bibitem{Burke2016Epi}
{\sc J.~V. Burke and T.~Hoheisel}, {\em Epi-convergence properties of smoothing
  by infimal convolution}, Set-Valued and Variational Analysis, 25 (2017),
  pp.~1--23.

\bibitem{cam}
{\sc L.~Cam}, {\em Th{\'e}orie asymptotique de la d{\'e}cision statistique},
  S{\'e}minaire de math{\'e}matiques sup{\'e}rieures, Presses de l'Univ. de
  Montr{\'e}al, 1972, \url{https://books.google.com/books?id=EH5EOwAACAAJ}.

\bibitem{dennismore}
{\sc J.~E. Dennis and J.~J. Mor{\'e}}, {\em A characterization of superlinear
  convergence and its application to quasi-newton methods}, Mathematics of
  computation, 28 (1974), pp.~549--560.

\bibitem{Eggermont2001}
{\sc P.~P.~B. Eggermont, V.~N. LaRiccia, and V.~LaRiccia}, {\em Maximum
  penalized likelihood estimation}, vol.~1, Springer, 2001.

\bibitem{fan1999one}
{\sc J.~Fan and J.~Chen}, {\em One-step local quasi-likelihood estimation},
  Journal of the Royal Statistical Society: Series B (Statistical Methodology),
  61 (1999), pp.~927--943.

\bibitem{ferguson}
{\sc T.~S. Ferguson}, {\em A course in large sample theory}, Routledge, 2017.

\bibitem{FriedlanderGoh}
{\sc M.~P. Friedlander and G.~Goh}, {\em Efficient evaluation of scaled
  proximal operators}, Electronic Transactions on Numerical Analysis, 46
  (2017), pp.~1--22.

\bibitem{glmnet}
{\sc J.~Friedman, T.~Hastie, and R.~Tibshirani}, {\em Regularization paths for
  generalized linear models via coordinate descent}, Journal of statistical
  software, 33 (2010), p.~1.

\bibitem{golub1966note}
{\sc G.~H. Golub and J.~H. Wilkinson}, {\em Note on the iterative refinement of
  least squares solution}, Numerische Mathematik, 9 (1966), pp.~139--148.

\bibitem{HareSagastizabal}
{\sc W.~Hare and C.~Sagastiz{\'a}bal}, {\em Computing proximal points of
  nonconvex functions}, Mathematical Programming, 116 (2009), pp.~221--258.

\bibitem{Hazan2016}
{\sc E.~Hazan, K.~Y. Levy, and S.~Shalev-Shwartz}, {\em On graduated
  optimization for stochastic non-convex problems}, in International conference
  on machine learning, 2016, pp.~1833--1841.

\bibitem{huang2019distributed}
{\sc C.~Huang and X.~Huo}, {\em A distributed one-step estimator}, Mathematical
  Programming, 174 (2019), pp.~41--76.

\bibitem{Ionides}
{\sc E.~Ionides}, {\em Maximum smoothed likelihood estimation}, Statistica
  Sinica,  (2005), pp.~1003--1014.

\bibitem{Kan21globalized}
{\sc C.~Kanzow and T.~Lechner}, {\em Globalized inexact proximal newton-type
  methods for nonconvex composite functions}, Computational Optimization and
  Applications, 78 (2021), pp.~377--410.

\bibitem{LeCam12}
{\sc L.~{Le Cam}}, {\em Asymptotic methods in statistical decision theory},
  Springer Science \& Business Media, 2012.

\bibitem{LeCamOG}
{\sc L.~LeCam}, {\em Locally asymptotically normal families of distributions},
  Univ. California Publ. Statist., 3 (1960), pp.~37--98.

\bibitem{LeCam70}
{\sc L.~LeCam}, {\em On the assumptions used to prove asymptotic normality of
  maximum likelihood estimates}, The Annals of Mathematical Statistics, 41
  (1970), pp.~802--828.

\bibitem{leesunsaunders}
{\sc J.~D. Lee, Y.~Sun, and M.~A. Saunders}, {\em Proximal newton-type methods
  for minimizing composite functions}, SIAM Journal on Optimization, 24 (2014),
  pp.~1420--1443.

\bibitem{eu-core}
{\sc J.~Leskovec, J.~Kleinberg, and C.~Faloutsos}, {\em Graph evolution:
  Densification and shrinking diameters}, ACM transactions on Knowledge
  Discovery from Data (TKDD), 1 (2007), pp.~2--es.

\bibitem{linearnonlinear}
{\sc D.~G. Luenberger, Y.~Ye, et~al.}, {\em Linear and nonlinear programming},
  vol.~2, Springer, 1984.

\bibitem{mattingley2012cvxgen}
{\sc J.~Mattingley and S.~Boyd}, {\em Cvxgen: A code generator for embedded
  convex optimization}, Optimization and Engineering, 13 (2012), pp.~1--27.

\bibitem{Mil16numerical}
{\sc A.~Milzarek}, {\em Numerical methods and second order theory for nonsmooth
  problems}, PhD thesis, Technische Universit{\"a}t M{\"u}nchen, 2016.

\bibitem{Mobahi}
{\sc H.~Mobahi and J.~W. Fisher}, {\em On the link between gaussian homotopy
  continuation and convex envelopes}, in International Workshop on Energy
  Minimization Methods in Computer Vision and Pattern Recognition, Springer,
  2015, pp.~43--56.

\bibitem{PaBo14proximal}
{\sc N.~Parikh, S.~Boyd, et~al.}, {\em Proximal algorithms}, Foundations and
  Trends{\textregistered} in Optimization, 1 (2014), pp.~127--239.

\bibitem{Pollard}
{\sc D.~Pollard}, {\em Another look at differentiability in quadratic mean}, in
  Festschrift for Lucien Le Cam, Springer, 1997, pp.~305--314.

\bibitem{polson2015proximal}
{\sc N.~G. Polson, J.~G. Scott, B.~T. Willard, et~al.}, {\em Proximal
  algorithms in statistics and machine learning}, Statistical Science, 30
  (2015), pp.~559--581.

\bibitem{VarAn}
{\sc R.~T. Rockafellar and R.~J.-B. Wets}, {\em Variational analysis},
  vol.~317, Springer Science \& Business Media, 2009.

\bibitem{lbfgs}
{\sc K.~Scheinberg and X.~Tang}, {\em Practical inexact proximal quasi-newton
  method with global complexity analysis}, Mathematical Programming, 160
  (2016), pp.~495--529.

\bibitem{Spokoiny}
{\sc V.~Spokoiny et~al.}, {\em Parametric estimation. finite sample theory},
  The Annals of Statistics, 40 (2012), pp.~2877--2909.

\bibitem{taddy2017one}
{\sc M.~Taddy}, {\em One-step estimator paths for concave regularization},
  Journal of Computational and Graphical Statistics, 26 (2017), pp.~525--536.

\bibitem{Tse09coordinate}
{\sc P.~Tseng and S.~Yun}, {\em A coordinate gradient descent method for
  nonsmooth separable minimization}, Mathematical Programming, 117 (2009),
  pp.~387--423.

\bibitem{vanderVaart}
{\sc A.~W. {Van der Vaart}}, {\em Asymptotic statistics}, vol.~3, Cambridge
  university press, 2000.

\bibitem{Xu2014Smoothing}
{\sc M.~Xu, J.~J. Ye, and L.~Zhang}, {\em Smoothing sqp methods for solving
  degenerate nonsmooth constrained optimization problems with applications to
  bilevel programs}, SIAM Journal on Optimization, 25 (2015), pp.~1388--1410.

\bibitem{zou2008}
{\sc H.~Zou and R.~Li}, {\em One-step sparse estimates in nonconcave penalized
  likelihood models}, Annals of statistics, 36 (2008), p.~1509.

\end{thebibliography}

\appendix
\section{Proofs}

\noindent 

\noindent\textit{Proof of Proposition \ref{stop_cond}.}
  Let $\Delta \hat{\theta}_{init} = \hat{\theta}_{ose} - \hat{\theta}_{init}$. It is a property of scaled proximal gradient descent that\footnote{See for example \cite{leesunsaunders}.}
  $$ C_{n} \Delta \hat{\theta}_{init} \in - \nabla g_{n}(\hat{\theta}_{init}) - \partial h_{n}(\hat{\theta}_{init} + \Delta \hat{\theta}_{init}).$$
So for some $u \in \partial h_{n}(\hat{\theta}_{init} + \nabla \hat{\theta}_{init})$ and $v \in \partial h_{n}(\hat{\theta})$
  $$\left(-\Delta \hat{\theta}_{init}\right)^{T} C_{n} \left( \hat{\theta}_{init} + \Delta \hat{\theta}_{init} - \hat{\theta}\right)$$
  $$= \left(\nabla g_{n}(\hat{\theta}_{init}) + u\right)^{T} \left(\hat{\theta}_{init} + \Delta \hat{\theta}_{init} - \hat{\theta}\right)$$
  $$= \left(\nabla g_{n}(\hat{\theta}_{init}) - \nabla g_{n}(\hat{\theta}) + u - v \right)^{T} \left(\hat{\theta}_{init} + \Delta \hat{\theta}_{init} - \hat{\theta}\right)$$
  By montonicity of the subdifferential of a convex function,
  $$ \geq \left(\nabla g_{n}(\hat{\theta}_{init}) - \nabla g_{n}(\hat{\theta})\right)^{T}\left(\hat{\theta}_{init} + \Delta \hat{\theta}_{init} - \hat{\theta}\right).$$
  By the strong convexity of $f_{n}$, with high probability and for some constant $m$ we can bound the previous display with
  $$ \geq \frac{m}{2} \left\| \hat{\theta}_{init} - \hat{\theta} \right\|^{2} + \left(\nabla g_{n}(\hat{\theta}_{init}) - \nabla g_{n}(\hat{\theta})\right)^{T} \Delta \hat{\theta}_{init}.$$
We have established
  $$(-\Delta \hat{\theta}_{init})^{T} \left(C_{n}(\hat{\theta}_{init} + \Delta \hat{\theta}_{init} - \hat{\theta}) - \left(\nabla g_{n}(\hat{\theta}_{init}) - \nabla g_{n}(\hat{\theta})\right) \right) \geq \frac{m}{2} \left\| \hat{\theta}_{init} - \hat{\theta} \right\|^{2}.$$
Thus
  $$- (\Delta \hat{\theta}_{init})^{T} C_{n} (\Delta \hat{\theta}_{init}) - (\Delta \hat{\theta}_{init})^{T} C_{n} (\hat{\theta}_{init} - \hat{\theta}) + (\Delta \hat{\theta}_{init})^{T} \left(\nabla g_{n}(\hat{\theta}_{init}) - \nabla g_{n}(\hat{\theta})\right)$$
  $$\geq \frac{m}{2} \left\| \hat{\theta}_{init} - \hat{\theta} \right\|^{2}.$$
Simplifying,
  $$-(\Delta \hat{\theta}_{init})^{T} C_{n} (\Delta \hat{\theta}_{init}) + (\Delta \hat{\theta}_{init})^{T} \left(\nabla g_{n}(\hat{\theta}_{init}) - \nabla g_{n}(\hat{\theta}) - C_{n}(\hat{\theta}_{init} - \hat{\theta})\right) \geq \frac{m}{2} \left\| \hat{\theta}_{init} - \hat{\theta} \right\|^{2}$$
so that Cauchy-Schwarz gives
  $$\|\Delta \hat{\theta}_{init} \|_{C_{n}}^{2} + \|\Delta \hat{\theta}_{init} \| \left\|\nabla g_{n}(\hat{\theta}_{init}) - \nabla g_{n}(\hat{\theta}) - C_{n}(\hat{\theta}_{init} - \hat{\theta})\right\| \geq \frac{m}{2} \left\| \hat{\theta}_{init} - \hat{\theta} \right\|^{2}.$$
  If $\lambda_{max}(C_{n}) < L$,
  $$L \| \Delta \hat{\theta}_{init} \|^{2} + \|\Delta \hat{\theta}_{init} \| \left\|\nabla g_{n}(\hat{\theta}_{init}) - \nabla g_{n}(\hat{\theta}) - C_{n}(\hat{\theta}_{init} - \hat{\theta})\right\| \geq \frac{m}{2} \left\| \hat{\theta}_{init} - \hat{\theta} \right\|^{2}.$$

  We have $\| \nabla g_{n}(\hat{\theta}_{init}) - \nabla g_{n}(\hat{\theta}) \| \leq M \| \hat{\theta}_{init} - \hat{\theta} \|$ by Lipschitz continuity of $\nabla g_{n}$. Thus
  $$L \| \Delta \hat{\theta}_{init} \|^{2} + M \|\Delta \hat{\theta}_{init} \| \cdot \|\hat{\theta}_{init}  - \hat{\theta} \|  + L \|\Delta \hat{\theta}_{init} \| \cdot \|\hat{\theta}_{init}  - \hat{\theta} \|$$
  $$\geq L \| \Delta \hat{\theta}_{init} \|^{2} + \|\Delta \hat{\theta}_{init} \| \left(M  \|\hat{\theta}_{init}  - \hat{\theta} \| + \| C_{n}(\hat{\theta}_{init} - \hat{\theta})\| \right)$$
  $$\geq L \| \Delta \hat{\theta}_{init} \|^{2} + \|\Delta \hat{\theta}_{init} \| \left( \|\nabla g_{n}(\hat{\theta}_{init}) - \nabla g_{n}(\hat{\theta}) \| + \| C_{n}(\hat{\theta}_{init} - \hat{\theta})\| \right)$$
  $$\geq L \| \Delta \hat{\theta}_{init} \|^{2} + \|\Delta \hat{\theta}_{init} \| \cdot \|\nabla g_{n}(\hat{\theta}_{init}) - \nabla g_{n}(\hat{\theta}) - C_{n}(\hat{\theta}_{init} - \hat{\theta})\| $$
  $$ \geq \frac{m}{2} \| \hat{\theta}_{init} - \hat{\theta} \|^{2}$$

Therefore
  $$L \| \Delta \hat{\theta}_{init} \|^{2} + M \|\Delta \hat{\theta}_{init} \| \cdot \|\hat{\theta}_{init}  - \hat{\theta} \|  + L \|\Delta \hat{\theta}_{init} \| \cdot \|\hat{\theta}_{init}  - \hat{\theta} \| \geq \frac{m}{2} \| \hat{\theta}_{init} - \hat{\theta} \|$$
  If $\|\Delta \hat{\theta}_{init}\| \geq \| \hat{\theta}_{init}- \hat{\theta} \|$ then
  $$ (2L + M)\| \Delta \hat{\theta}_{init} \|^{2} \geq \frac{m}{2} \| \hat{\theta}_{init} - \hat{\theta} \|^{2}$$
  $$\Rightarrow \| \Delta \hat{\theta}_{init} \| \geq \sqrt{\frac{m}{2(2L + M})} \| \hat{\theta}_{init} - \hat{\theta} \|.$$
  Otherwise, if $\|\Delta \hat{\theta}_{init}\| < \| \hat{\theta}_{init}- \hat{\theta} \|$ then
  $$(2L + M) \|\Delta \hat{\theta}_{init} \| \cdot \| \hat{\theta}_{init} - \hat{\theta} \| \geq \frac{m}{2} \| \hat{\theta}_{init} - \hat{\theta} \|^{2}$$
  $$\Rightarrow \|\Delta \hat{\theta}_{init} \| \geq \frac{m}{2(2L + M)} \| \hat{\theta}_{init} - \hat{\theta} \|$$

  Therefore $\|\Delta \hat{\theta}_{init}\| \geq \max\left\{ \sqrt{\frac{m}{2(2L + M)}}, \frac{m}{2(2L + M)} \right\} \|\hat{\theta}_{init} - \hat{\theta} \|$ is small.

  Recalling that $\Delta \hat{\theta}_{init} = \hat{\theta}_{ose} - \hat{\theta}_{init}$, the result is proven. 
\qed
\endproof

  \smallskip

\vspace{.5cm}

\noindent\textit{Proof of Proposition~\ref{prop:diffeCf}.}

Let $x^*\in\R^d$, and let $w^*\in \prox_{f,C}(x^*)$ be defined as in the statement for $u^\nu\to 0$, and $w^\nu\in\prox_{f,C}(x^*+u^\nu)$, with $w^\nu \to w^*$.  
We need to show that
\[\lim_{u \to 0} \frac{e_C f(x^*+u)-e_C f(x^*)-u^\top C(x^*-w^*)}{\|u\|}=0\]

By definition of the Moreau envelope,
\begin{align*}
  e_C f(x^*+u^\nu)&=f(w^\nu)+\frac{1}{2}(w^\nu-(x^*+u^\nu))^\top C(w^\nu-(x^*+u^\nu))\\
  & \leq f(w^*)+\frac{1}{2}(w^*-(x^*+u^\nu))^\top C(w^*-(x^*+u^\nu))
\end{align*}
Hence
\begin{align*}
  &\quad \frac{e_C f(x^*+u^\nu)-e_C f(x^*)-(u^{\nu})^\top C(x^*-w^*)}{\|u^\nu\|}\\
  &\leq \frac{f(w^*)+\frac{1}{2}(w^*-(x^*+u^\nu))^\top C(w^*-(x^*+u^\nu))-e_{C} f(x^{*})-(u^{\nu})^\top C(x^{*}-w^{*})}{\|u^{\nu}\|}\\
  &\leq \frac{\frac{1}{2}(u^\nu)^\top C(w^*-x^*)-(u^{\nu})^{\top} C(x^*-w^*)}{\|u^{\nu}\|}\\
& = \frac{(u^\nu)^\top C u^\nu}{\|u^{\nu}\|}
\end{align*}
where the second inequality follows from the definition  of $e_C f(x^*)=f(w^*)+\frac{1}{2}(w^*-x^*)^\top C(w^*-x^*)$.

On the other hand, a lower bound on the quotient follows from
$$e_C f(x^*)=f(w^*)+\frac{1}{2}(w^*-x^*)^\top C(w^*-x^*)
\leq f(w^\nu)+\frac{1}{2}(w^\nu-x^*)^\top C(w^\nu-x^*)$$
Hence,
\begin{align*}
  &\frac{e_C f(x^*+u^\nu)-e_C f(x^*) -(u^{\nu})^{\top} C(x^*-w^*)}{\|u^\nu\|}\\
  &\geq \frac{e_C f(x^*+u^\nu)-\left(f(w^\nu)+\frac{1}{2}(w^\nu-x^*)^\top C(w^\nu-x^*)\right) -(u^{\nu})^{\top} C(x^*-w^*)}{\|u^\nu\|}\\
&=\frac{\frac{1}{2}(u^\nu)^\top Cu^\nu-(u^\nu)^\top C(w^\nu-x^*)-(u^\nu)^\top C(x^*-w^*)}{\|u^\nu\|}\\
&=\frac{\frac{1}{2}(u^\nu)^\top Cu^\nu+(u^\nu)^\top C(w^*-w^\nu)}{\|u^\nu\|}
\end{align*}
Combining both bounds, we get the following inequality chain.
$$\frac{\frac{1}{2}(u^\nu)^\top C u^\nu+(u^\nu)^\top C (w^*-w^\nu)}{\|u^\nu\|}$$
$$\leq \frac{e_C f(x^*+u^\nu)-e_C f(x^*) -(u^{\nu})^{\top} C(x^*-w^*)}{\|u^\nu\|} $$
$$\leq \frac{(u^\nu)^\top C u^\nu}{\|u^{\nu}\|}.$$
Taking the limit when $u^\nu\to 0$, the lower bound goes to zero as $w^\nu\to w^*$, whereas the upper bound also goes to zero as $\nu \to \infty$. Since this inequality holds for any selection of $u^\nu\to 0$, this establishes that $e_C f$ is differentiable, with gradient $\nabla e_C f(x^*) = C(x^*-w^*)$.
%(\cite[Prop.12.30]{BauCom11cvx})
\qed
\endproof

\end{document}